\documentclass{amsart}

\usepackage{amsthm}

\newtheorem{lemma}{Lemma}
\newtheorem{proposition}{Proposition}
\newtheorem{theorem}{Theorem}
\newtheorem{corollary}{Corollary}

\newtheorem{remark}{Remark}

\newcommand{\nth}[1]{$#1 {\rm - th }$}

\newcommand{\Z}{\mathbb{Z}}
\newcommand{\ZM}[1]{\Z /( #1 \cdot \Z)}
\newcommand{\ZMs}[1]{(\Z / #1 \cdot \Z)^*}

\newcommand{\rk}[1]{{\rm rank}(#1)}
\newcommand{\ord}{{\rm ord}}

\newcommand{\Tr}{\mbox{\bf Tr}}

\newcommand{\rg}[1]{\mbox{\bf #1}}
\newcommand{\eu}[1]{\mathfrak{#1}}
\newcommand{\id}[1]{\mathcal{#1}}
\newcommand{\Gal}{\mbox{ Gal }}
\newcommand{\diag}{\mbox{ \bf{Diag} }}
\newcommand{\rf}[1]{(\ref{#1})}

\newcommand{\nmid}{\not \hspace{0.25em} \mid}
\newcommand{\Norm}{\mbox{\bf N}}

\newcommand{\F}{\mathbb{F}}
\newcommand{\K}{\mathbb{K}}
\newcommand{\LK}{\mathbb{L}}
\newcommand{\KL}{\mathbb{L}}

\newcommand{\Q}{\mathbb{Q}}
\newcommand{\R}{\mathbb{R}}
\newcommand{\C}{\mathbb{C}}
\newcommand{\N}{\mathbb{N}}

\begin{document}
{\obeylines \small
\vspace*{0.2cm}
\hspace*{6.0cm}Ils ne savent pas ce qu'ils perdent
\hspace*{6.0cm}Tous ces sacr\'es cabotins.
\hspace*{6.0cm}Sans le latin, sans le latin,
\hspace*{6.0cm}La messe nous emmerde. \footnote{Georges Brassens}
\vspace*{0.8cm}
\hspace*{7.5cm}\indent{\it To Seraina and Theres}
\vspace*{1.5cm}
\smallskip}
\title[Catalan-Fermat]
{A Cyclotomic Investigation of the Catalan -- Fermat Conjecture. DRAFT}
\author{Preda Mih\u{a}ilescu}
\address[P. Mih\u{a}ilescu]{Gesamthochschule Paderborn}
\email[P. Mih\u{a}ilescu]{preda@upb.de}
\date{Version 1.0 \today}
\vspace{2.0cm}
%\huge
%%%%%%%%%%%%%%%%%%%%%%%%%%%%%%%%%%%%%%%%%%%%%%%%%%%%%%%%%%%%%%%%%%
\begin{abstract}
With give some new, simple results on the equation $x^{p}+y^{p} = z^{q}$
using classical methods of cyclotomy.
\end{abstract}
\tableofcontents
\maketitle
\vspace*{1.0cm}
\section{Introduction}
Consider the equation 
\begin{eqnarray}
\label{FC}
x^{p} + y^{p} = z^{q}, \quad \hbox{with} \quad  x, y, z \in \Z, \ (x,
y, z) = 1 \quad \hbox{and $p, q$ odd primes.}
\end{eqnarray}
This is the natural generalization the equations $x^{p} + y^{p} +
z^{p} = 0$ of Fermat and $x^{p} - y^{q} = 1$ of Catalan, and we shall
denote it by \textit{Fermat - Catalan} equation. This name is used by
some authors for the more general
\begin{eqnarray}
\label{SF}
x^{p} + y^{q} = z^{r},
\end{eqnarray} 
while others \cite{Za} use the term of \textit{super - Fermat}
equation, for \rf{SF}. We shall also adopt this terminology in this
paper. There has been an increasing literature on the subject in the
last decade, yet general results are still scarce. Thus F. Buekers
enumerates all the solutions of \rf{SF} for $\chi = 1/p + 1/q + 1/r
\geq 1$ , while for $\chi < 1$, Darmon and Granville \cite{DG} have
proved that there are at most finitely many solutions for a fixed set of
exponents. Specific curves are used for fixed triples of exponents by
Darmon and coauthors \cite{Da1}, \cite{DM}, Ellenberg \cite{El}, Bruin
\cite{Br}, Poonen et. al. \cite{PSS}. We refer the reader to \cite{Be}
for a nice survey of the topic and a comprehensive overview of the
cases known up to fall 2004.

A program generalizing the Wiles proof of Fermat's Last Theorem was
proposed by Darmon in \cite{Da}; it suggests replacing elliptic curves
by hyperelliptic curves or even surfaces, in the general case.  A
particular case we shall consider is the \textit{rational} Catalan
equation:
\begin{eqnarray}
\label{ratcat}
X^{p} + Y^{q} = 1, \quad \hbox{with odd primes $p, q$ and } \quad X, Y
\in \Q; 
\end{eqnarray}
it is easily shown (see below) that this is equivalent to the special
case $x^{p}+y^{q}= z^{pq}$ of \rf{FC}.  Tijdeman and Shorey emitted in
\cite{ST}, Chapter XII, the conjecture, that \rf{ratcat} has at most
finitly many rational solutions. We shall find conditions on $(p, q)$
for which that equation has no non-trivial rational solutions at all.

Our purpose in this paper is to investigate \rf{FC} separately, using
classical cyclotomic approaches. Unsurprisingly, the results obtained
are partial, but the set of exponents $p, q$ for which they are valid
are unbounded. The results are ordered in increasing order of
conditions on the (odd prime) exponents $p, q$. Based on a simple
relative class number divisibility condition we reduce first \rf{FC}
to a Fermat - type equation (with only one exponent) over a totally
real field:

\begin{theorem}
\label{genfer}
Let $p, q$ be odd primes with $p >3, q \nmid h_{p}$ - with $h_{p}$ the
class number of the \nth{p} cyclotomic extension - and for which the
Fermat - Catalan equation \rf{FC} has a non trivial solution. Then the
equation
\begin{eqnarray}
\label{fereq}
a X^{q} + b Y^{q} + c Z^{q}  =  0 \quad & & X \in \Z; \quad  Y, Z \in
\Z[\zeta+\overline \zeta] 
\end{eqnarray}
has a non - trivial solution. Here $b, c$ are units, while $a$ is
either unit or the principal ideal $(a)$ is a power of the ramified
prime ideal above $p$.
\end{theorem} 
After deriving the distinction of the first and second case in
\rf{FC}, and the analogs of the Barlow - Abel formulae \cite{Ri}, we
prove a theorem, which allows a useful additional case distinction:

\begin{theorem}
\label{main}
If the equation \rf{FC} has a solution $(x,y,z; p, q)$ and the
exponents are such that $\max\{p,\frac{p(p-20)}{16}\} > q$ and $q \nmid
h_{p}^{-}$, with $h_{p}^{-}$ the relative class number of the \nth{p}
cyclotomic extension, then
\begin{eqnarray}
\label{qpot2}
x+f \cdot y & \equiv & 0 \mod q^{2} \quad \hbox{ with } \quad f \in
\{-1,0,1\}. 
\end{eqnarray}
\end{theorem}

This leads to a delicate case by case analysis (six cases in
total). For five of the six cases we are able to find some simple
algebraic conditions, while one of the cases ($z \not \equiv 0 \mod p$
and $x \equiv 0 \mod q^{2}$) remains unsolved. The main results of
this paper are the following:

\begin{theorem}
\label{main1}
Let $p, q$ be odd primes such that $q \nmid h(p,q); \ p \not \equiv 1
\mod q$ and $\max\{p,\frac{p(p-20)}{16}\} > q$. Then the equation
\[ x^{p} + y^{p} = z^{q} , \quad (x, y, z) = 1 \]
has no non - trivial solutions for 
\begin{eqnarray}
\label{lowb}
\max\{|x|,|y|\} \geq \frac{1}{2} \cdot \left(\frac{1}{p(p-1)} \cdot
\left(\frac{q^{\frac{q-2}{q-1}}}{2 }\right)^{p-2}\right)^{q}.
\end{eqnarray}
 \end{theorem}
 Here $h(p,q)$ is a divisor of the relative class number $h_{pq}^{-}$
 which is explicitly computable by means of generalized Bernoulli
 numbers and will be defined below. Assuming the stronger condition
 $q \nmid h_{pq}^{-}$, the relative class number of the \nth{pq}
 cyclotomic field, we have:
\begin{theorem}
\label{main2}
Let $p, q > 3$ be primes such that \rf{FC} has a solution and suppose
that $-1 \in < p \mod q>$, $\max\{p,\frac{p(p-20)}{16}\} > q$ and
$q \nmid h_{pq}^-$. Then either
\begin{eqnarray}
\label{Wifcases}
a^{q-1} \equiv 1 \mod q^{2} \quad \hbox{ for some } \quad a \in \{ 2, p,
2^{p-1} \cdot p^{p} \},
\end{eqnarray}
or 
\begin{itemize}
\item[A.] $p \nmid z$ and $q^{2} | x y$ if $q \not \equiv 1 \mod p$ and
  $q^{3} | x y$, if $q \equiv 1 \mod p$.
\item[B.] If $q \not \equiv 1 \mod p$, then 
\begin{eqnarray}
\label{bound}
\min(|x|, |y|) > c_{1}(q) \left(\frac{q^{p-1}}{p} \right)^{q-2}, \quad
\hbox{if} \quad q \not \equiv 1 \mod p,
\end{eqnarray}
and
\begin{eqnarray}
\label{bound1}
\min(|x|, |y|) > c_{1}(q) \left(\frac{q^{2(p-1)}}{p} \right)^{q-2},
\end{eqnarray}
otherwise. Here $c_{1}(q)$ is am effectively computable, strictly
increasing function with $c_{1}(5) > 1/2$.
\end{itemize}
\end{theorem}

An immediate consequence of this Theorem is the following
generalization of Catalan's conjecture:
\begin{corollary}
\label{catg1}
Let $p, q$ be odd primes such that 
\begin{enumerate}
\item $-1 \not \in <p \mod q> $ and $q \nmid h_{pq}^-$,
\item  $\max\{p,\frac{p(p-20)}{16}\} > q$, 
\item $a^{q-1} \not \equiv 1 \mod q^{2}$ for $a \in \{ 2, p, 2^{p-1}
  \cdot p^{p} \}$.
\end{enumerate}
 Then the equation
\[ X^{p} + C^{p} = Z^{q} \]
has no integer solution for fixed $C$ with $|C| < \frac{1}{2} \cdot
\left(\frac{q^{p-1}}{p} \right)^{q-2}$. If $q \equiv 1 \mod p$, there
are no solutions with $|C| < \frac{1}{2} \cdot
\left(\frac{q^{2(p-1)}}{p} \right)^{q-2}$
\end{corollary}
\begin{proof}
The premises allow us to apply Theorem \ref{main2} and the claim
follows from \rf{bound} and \rf{bound1}.
\end{proof}

Finally, by restricting the results above to the rational Catalan
equation \rf{ratcat}, we are able to give a criterion which only
depends on the exponents, namely the following:
\begin{theorem}
\label{trc}
Let $p, q > 3$ be distinct primes for which the following conditions
are true:
\begin{itemize}
\item[1.] $-1 \in < p \mod q>$ and $-1 \in < q  \mod p>$,
\item[2.] $\left(pq,  h_{pq}^-\right) = 1$,
\item[3.] $2^{p-1} \not \equiv 1 \mod p^{2}$ and $2^{q-1} \not \equiv
1 \mod q^{2}$, 
\item[4.] $\left(2^{p-1} p^{p}\right)^{q-1} \not \equiv 1 \mod q^{2}$
  and $\left(2^{q-1} q^{q}\right)^{p-1} \not \equiv 1 \mod p^{2}$,
\item[5.] $p^{q-1} \not \equiv 1 \mod q^{2}$ and $q^{p-1}
  \not \equiv 1 \mod p^{2}$, 
\item[6.] $\max\{p,\frac{p(p-20)}{16}\} > q$ and
  $\max\{q,\frac{q(q-20)}{16}\} > p$. 
\end{itemize} 
Then the equation $X^{p}+Y^{q} = 1$ has no rational solutions.
\end{theorem}

\section{Generalities}
The following lemma of Euler is used in both equations of Fermat and
Catalan:
\begin{lemma}
\label{euler}
Let $x, y$ be coprime integers and $n > 1$ be odd. Then 
\begin{eqnarray}
\label{cop}
\left(\frac{x^{n}+y^{n}}{x+y}, x+y\right) \ | \ n.
\end{eqnarray}
\end{lemma}
\begin{proof}
Write $x = (x+y) - y$ and develop
\begin{eqnarray}
\label{dev} 
\frac{x^{n}+y^{n}}{x+y} & = & \frac{(x+y)^{n} + \sum_{k=1}^{n-2}
  \binom{n}{k} (x+y)^{n-k} \cdot (-y)^{k}  + n (x+y) y^{n-1} - y^{n} +
  y^{n}}{x+y} \nonumber \\
& = & K \cdot (x+y) + n \cdot y^{n-1}, \quad \hbox{with } \quad K \in \Z.
\end{eqnarray}
The common divisor in \rf{cop} is consequently 
\[ D = \left(\frac{x^{n}+y^{n}}{x+y}, x+y\right) = \left(K \cdot (x+y) + n
  \cdot y^{n-1}, x+y\right) = \left(n \cdot y^{n-1}, x+y\right). \] But
  since $x, y$ are coprime, and consequently also $(x+y, y) = 1$, it
  follows plainly that $D = (n, x+y) | n$.
\end{proof}
Like in Fermat's equation, the above lemma leads to the following case
distinction for \rf{FC}:

\textbf{Case I}: The case in which $p \nmid z$. Then 
\[ x^{p} + y^{p} = (x+y) \cdot \frac{x^{p}+y^{p}}{x+y} = z^{q}, \]
and since by Lemma \ref{euler}, the two factors have no common divisor,
they must simultaneously be \nth{q} powers. Thus
\begin{eqnarray}
\label{case1}
x+y = A^{q}, \quad \frac{x^{p}+y^{p}}{x+y} = B^{q} \quad \hbox{with}
\quad A, B \in \Z,
\end{eqnarray}
for this case.

\textbf{Case II}: The case in which $p \mid z$. Then $x^{p} + y^{p}
\equiv x+y \equiv z^{q} \equiv 0 \mod p$. We show that in this case 
\[ v_{p}\left( \frac{x^{p}+y^{p}}{x+y}\right)   = 1. \]
The development in \rf{dev} yields
\[ \frac{x^{p}+y^{p}}{x+y} = p \cdot y^{p-1} + \sum_{k=2}^{p-1}
\binom{p}{k} (x+y)^{k-1} \cdot (-y)^{p-k} + (x+y)^{p-1}. \] The
binomial coefficients in the above sum are all divisible by $p$. Since
$(y, x) = (y, x+y) = 1$, it follows also that $(y,p) = 1$ and thus
$v_{p}(p y^{p-1}) = 1$. But $p | x+y$, so all the remaining terms in
the expansion of $\frac{x^{p}+y^{p}}{x+y}$ are divisible (at least) by
$p^{2}$, which confirms our claim. In the second case thus, a \nth{q}
power of $p$ is split between the factors $x+y$ and
$\frac{x^{p}+y^{p}}{x+y}$ in such a way that the latter is divisible
exactly by $p$. In the second case we have herewith: 
\begin{eqnarray}
\label{case2}
x+y = p^{nq-1} \cdot A^{q} = (A')^{q}/p^{e}, & \quad &
\frac{x^{p}+y^{p}}{x+y} = p \cdot B^{q}  \\
 \hbox{with}  \quad (A, B, p) = 1 & \quad & \hbox{and}  \quad n =
 v_{p}(z) \geq 1. \nonumber
\end{eqnarray}

The relations \rf{case1} and \rf{case2} are the analogs of the Barlow
- Abel relations for the Fermat - Catalan equations. Fermat's equation
is homogeneous; thus if $x^{p}+y^{p}+z^{p} = 0$ is a solution in which
$(x,y,z)$ are not coprime, one may divide by the \nth{p} power of the
common divisor, thus obtaining a solution with coprime $(x,y,z)$. The
equation \rf{FC} is {\em not} homogeneous and thus one may ask whether
the requirement that $(x,y,z)$ be coprime is not
restrictive\footnote{I owe David Masser the observation that the
  condition $(x,y,z) = 1$ is not obvious for the equation \rf{FC}}.
The following lemma addresses this question and shows that one can
construct arbitrary many solutions of $x^{p}+y^{p} = z^{q}$, if common
divisors are allowed. It appears that the condition $(x,y,z) = 1$ is
thus plausible.

\begin{lemma}
  Let $x, y$ be coprime integers and $p, q$ be distinct odd primes.
  Then there is an integer $D \in \Z$ such that
\[ (D \cdot x)^{p} + (D \cdot y)^{p} = z^{q}, \quad \hbox{with} \quad
z \in \Z. \] Furthermore, every integer solution of
$X^{p}+Y^{p}=Z^{q}$ with $(X, Y, Z) > 1$ arises in this way.
\end{lemma}
\begin{proof}
Let $x, y$ be coprime integers and 
\[ x^{p} + y^{p} = C \cdot z^{q}, \]
where $C \in \Z$ and $z$ is the largest \nth{q} power dividing the
left hand side ($z = 1$ is possible); i.e. $C$ is \nth{q} power -
free. Let $\ell | C$ be a prime and $n = v_{\ell}(C)$. We show that
there is an integer $D(\ell)$ such that $x' = D(\ell) \cdot x$ and $y'
= D(\ell) \cdot y$ verify ${x'}^{p} + {y'}^{p} = C' \cdot {z'}^{q}$
and $z | z', C' | C$ while $(C', \ell) = 1$. Indeed, let $a \in \N$ be
such that $n + a\cdot p \equiv 0 \mod q$; such an integer exists since
$(p, q) = 1$.  We define $D(\ell) = \ell^{a}$ and $b = (n+ap)/q$. Then
\begin{eqnarray*}
{x'}^{p} + {y'}^{p} =  \ell^{ap} \cdot C \cdot z^{q} = \ell^{ap+n}
\cdot C/\ell^{n} \cdot z^{q} 
 =  \frac{C}{\ell^{a}} \cdot (z \cdot \ell^{b})^{q}.
\end{eqnarray*}
Setting $C' = C/\ell^{a}$ and $z' = z \cdot \ell^{b}$, the claim
follows. By repeating the procedure recursively for all the prime
divisors of $C$ one obtains $D = \prod_{\ell | C} \ D(\ell)$ for which
the claim of the Lemma holds.

Conversely, let $x^{p}+y^{p} = z^{q}$ hold for a triple with $(x,y,z)
= G$. For each prime $\ell | G$ let $a = v_{\ell}(x,y)$. Then $a p
\leq q \cdot v_{\ell}(z)$. If $G' = (x, y)$ it follows that $w =
z^{q}/{G'}^{p}$ is an integer. The integers $x' = x/G', y' = y/G'$ are
coprime; if $C$ is the \nth{q} power - free part of $w$ and ${z'}^{q}
= w/C$, then
\begin{eqnarray}
\label{above}
  {x'}^{p} + {y'}^{p} = C \cdot {z'}^{q} .
\end{eqnarray}
But then the initial equality can be derived from \rf{above} by the
procedure described above. Thus all non trivial solutions of
$x^{p}+y^{p} = z^{q}$ have coprime $x, y, z$.
\end{proof}
We finally prove the relation between \rf{ratcat} and \rf{FC}:
\begin{lemma}
\label{rcfc}
Let $p, q$ be odd primes. The equation \rf{ratcat} has non trivial
rational solutions if and only if
\begin{eqnarray}
\label{homfc}
x^{p}+y^{q} = z^{pq}, \quad (x, y, z) = 1 \quad \hbox{ and } \quad x,
y, z \in \Z
\end{eqnarray}
has non trivial solutions.
\end{lemma}
\begin{proof}
  Suppose first that \rf{homfc} has some non trivial solution $x, y,
  z$. Then one easily verifies that \rf{ratcat} has the solution $X =
  x/z^{q}, Y = y/z^{p}$. 

Conversely, let $\ X = a/c, Y = b/d$ be a non trivial solution of
\rf{ratcat} with $(a, c) = (b, d) = 1; \ a, b, c, d, \in
\Z$. Clearing denominators we find
\[ a^{p} d^{q} + b^{q} c^{p} = c^{p} d^{q}. \] 
Since $(a, c) = (b, d) = 1 $, by comparing the two sides of the identity,
we find $c^{p} | d^{q}$ and $d^{q} | c^{p}$; thus $c^{p} = d^{q}$. But
$p$ and $q$ are distinct primes, thus for each prime $\ell | c$, we
have $pq | v_{\ell}(c)$ and we may write $c^{p} = d^{q} = u^{pq}$. The
equation \rf{ratcat} becomes $a^{p}+b^{q}=u^{pq}$, as claimed.
\end{proof}
\section{Cyclotomy and Fermat - Catalan}
We start by fixing some notations which shall be used throughout the
rest of this paper.
\subsection{Notation}
We shall let $p, q$ be two odd primes and $\zeta, \xi \in \C$ be
primitive \nth{p} and \nth{q} roots of unity and $\K = \Q(\zeta), \K'
= \Q(\xi), \KL= \Q(\zeta, \xi)$ the respective cyclotomic fields.
Furthermore, the Galois groups will be
\begin{eqnarray*}
 G_{p} & = & \Gal(\K/\Q) = < \sigma >  ,  \quad G_{q} = \Gal(\K')/\Q = < \tau
 >\quad \hbox{ and } \\
 G & = & \Gal(\KL/\Q) = < \sigma \tau > = < \sigma > \times < \tau > . 
\end{eqnarray*} 
Unless stated otherwise in some particular context, $\sigma, \tau$ are
thus generators of the Galois groups $G_{p}, G_{q}$, respectively.
Furthermore, if $0< a < p; 0 < b < q$ we shall use the notation
$\sigma_{a} \in G_{p}, \tau_{b} \in G_{q}$ for the elements of the
Galois groups given by $\zeta \mapsto \zeta^{a}$ and $ \xi \mapsto
\xi^{b}$, respectively. The map $G_{p} \rightarrow \ZMs{p}$ given by
$\sigma_{a} \mapsto a$ will be denoted by $\widehat{\sigma_{a}} = a$,
and likewise for the analog for $G_{q}$. Complex conjugation is denoted
also by $\jmath \in G$, while $\jmath_{p}, \jmath_{q}$ are the complex
conjugation maps of $G_{p}, G_{q}$ respectively, lifted to $G$. Thus
$\jmath_{p}$ acts on $\zeta$ but fixes $\xi$ and $\jmath_{q}$ does the
reverse.  

Finally, the ramified primes are $\wp = (1-\zeta), \eu{q} = (1-\xi)$.
We use the notations $\lambda = (\xi - \overline \xi)$ and $\lambda' =
(\zeta - \overline \zeta)$ for generators of these ramified primes;
this is due to the nice behavior under complex conjugation. In several
contexts it will come natural to use the classical $\lambda =
(1-\xi)$, etc.; this deviation from the general use will be mentioned
in place.

We now start with some classical results, adapted to the present
equation \rf{FC}.
\subsection{First Consequences of Class Field Theory}
We assume $q \nmid h_{p}^{-}$ and deduce some consequences, starting
from a presumed non - trivial solution $(x, y, z)$ of \rf{FC}. With
$e$ defined like above, we shall let  
\begin{eqnarray}
\label{alphadef}
 \alpha & = & \frac{x+\zeta \cdot y}{(1-\zeta)^{e}} \quad \hbox{ and }
 \\
\eu{A} & = & \left(\alpha, \frac{x^{p}+y^{p}}{p^{e}(x+y)}\right)
 \subset \id{O}(\K)^{\times} = \Z[\zeta]^{\times}. \nonumber
\end{eqnarray}

\begin{lemma}
\label{euaq}
Let $x, y, z$ be coprime integers verifying \rf{FC}. Then 
\begin{eqnarray}
\label{firstalpha}
\left(\sigma(\alpha), \sigma'(\alpha)\right) &= & 1 \quad \forall \
\sigma, \sigma' \in G_{p}, \sigma \neq \sigma' \\
\eu{A}^{q} & = & (\alpha).\nonumber
\end{eqnarray}
\end{lemma}
\begin{proof}
We cumulate the two Cases of the Fermat - Catalan equation in
\begin{eqnarray}
\label{cum1}
 x+y & = & A^{q}/p^{e} \quad \hbox{and} \\
\label{cum2}
 \frac{x^{p}+y^{p}}{x+y} & = & p^{e} \cdot B^{q} \quad \hbox{with} \quad e
 \in \{ 0, 1\}
\end{eqnarray}
Here $e = 0$ corresponds to the First Case and $e = 1$ to the Second
Case. If $e = 1$, it is understood that $p | A$, so the right hand
side in \rf{cum1} is an integer. Then \rf{cum2} is equivalent to
\begin{eqnarray}
\label{neq}
\Norm_{\Q(\zeta)/\Q}(\alpha) = B^{q} \quad \hbox{ and } \quad \eu{A} =
(\alpha, B).
\end{eqnarray}
We write $P = \{1, 2, \ldots, p-1\}$ and note that, for $c, d \in P$
we have
\begin{eqnarray}
\label{copr}
\Delta(c,d) = (\sigma_{c}(\alpha), \sigma_{d}(\alpha)) = (1)
\end{eqnarray}
whenever $c \neq d$. Indeed 
\begin{eqnarray*}
 y(\zeta^{c}-\zeta^{d}) =
(1-\zeta^{c})^{e} \sigma_{c}(\alpha) - (1-\zeta^{d})^{e}
\sigma_{d}(\alpha) \in \Delta(c,d) \quad \hbox{and} \\
 x(\zeta^{-c}-\zeta^{-d}) =
\overline \zeta^{c} (1-\zeta^{c})^{e} \sigma_{c}(\alpha) - \overline
\zeta^{d }(1-\zeta^{d})^{e} \sigma_{d}(\alpha) \in \Delta(c,d).
\end{eqnarray*}
Since $(x, y) = 1$ it follows that $\Delta(c,d) \supset \wp$, the
ramified prime above $p$ in $\K$. But $\alpha = (x+\zeta
y)/(1-\zeta)^{e} = \frac{x+y}{(1-\zeta)^{e}} - y \cdot
(1-\zeta)^{1-e}$. If $e = 0$, the first term is coprime to $\wp$ and
the second is not. If $e = 1$, the second term is coprime to $\wp$ and
the first is not. Thus in both cases, $(\alpha, \wp) = 1$ and
$\Delta(c,d) = (1)$, as claimed. The second relation in
\rf{firstalpha} follows now easily from the definition of $\eu{A}$:
\[ \eu{A}^{q}/(\alpha) = \left(\alpha^{q-1}, B \cdot \alpha^{q-2},
  \ldots, B^{q-1}, \Norm(\alpha)/\alpha \right) , \] and one verifies
that the integer ideal on the right hand side is equal to the ideal
$\left(\alpha, \Norm(\alpha)/\alpha\right) = (1)$.
\end{proof}

An immediate consequence is:
\begin{corollary}
\label{crhodef}
If $p, q$ are odd primes with $q \nmid h_p$, the class number of
$\Q(\zeta_p)$, $x,y,z$ verify \rf{FC} and $\alpha = (x+\zeta
y)/(1-\zeta)^e$ with $e$ as above, then
\begin{eqnarray}
\label{algq}
\alpha & = & \varepsilon \cdot \rho^q \quad \hbox{for some} \quad
\varepsilon \in \Z[\zeta + \overline \zeta]^{\times} , \ \rho \in
\Z[\zeta].
\end{eqnarray}
and if only $q \nmid h_{p}^{-}$ holds, then
\begin{eqnarray}
\label{rhod}
\frac{x+\zeta y}{x + \overline \zeta y} = \pm \left(\frac{\rho_{1}}{\overline
    \rho_{1}}\right)^{q},  \quad \hbox{for some} \quad \rho_{1}
\Z[\zeta].
\end{eqnarray}
The two algebraic integers $\rho, \rho_{1}$ may, but need not be equal.
\end{corollary}
\begin{proof}
  The statement \rf{algq} is a direct consequence of \rf{firstalpha},
  since the second relation implies that $\eu{A}$ is principal if $q
  \nmid h_{p}$. Since $\eu{A}^{q} = (\alpha)$ by \rf{firstalpha}, if
  $q \nmid h_{p}^{-}$, then there is a real ideal $\eu{B} \subset
  \Z[\zeta]$ together with an algebraic number $\nu \in \K$ such that
  $\eu{A} = (\nu) \cdot \eu{B}$. By dividing through the complex
  conjugate of this identity, one finds
\[  \left(\frac{\eu{A}}{\overline{\eu{A}}}\right)^{q} =
(\alpha/\overline \alpha) = (\nu/\overline \nu)^{q}, \] and there is a
unit $\eta$ such that $\alpha/\overline \alpha = \eta (\nu/\overline
\nu)^{q}$ and thus
\[ \frac{x+\zeta y}{x+\overline \zeta y} = \left(\frac{1-\overline
    \zeta}{1-\zeta}\right)^{e} \cdot \eta \cdot
\left(\frac{\nu}{\overline \nu}\right)^{q} = \eta' \cdot
\left(\frac{\nu}{\overline \nu}\right)^{q}. \] But then $\eta' \cdot
\overline \eta' = 1$ and Dedekind's unit Theorem implies that $\eta'$
is a root of unity of $\K$. Since all roots of unity of this field
have order dividing $2p$, and $(2p, q) = 1$, the statement \rf{rhod}
follows.
\end{proof}

We now prove a lemma concerning ideals related to the above $\eu{A}$,
in a more general setting.
\begin{lemma}
\label{idealq}
Let $\rg{k} \subset \LK$ be some field such that $q \nmid h(\rg{k}$,
the class number of the field $\rg{k}$. Let $\phi_i \in
\id{O}(\rg{k}), i=1, 2, \ldots, n$ be such that $(\phi_i, \phi_j) =
(1)$ for $1 \leq i \neq j < n$.  Suppose that for some $m \geq 0$ all
$\phi_i, i = m+1, m+2, \ldots, $ are not units, while $\phi_j, j = 1,
2, \ldots, m$ are units.  Furthermore, there is a $C \in \rg{k}, (C,
pq) = 1$ such that
\begin{eqnarray}
\label{prd}
\prod_{i=1}^n \phi_i = C^{q}.
\end{eqnarray}
Then there are $\eta_i \in \id{O}(\rg{k})^{\times}$ and $\mu_i \in
\id{O}(\rg{k})$ such that
\begin{eqnarray}
\label{etamu}
\phi_i = \eta_i \cdot \mu_i^q \quad \hbox{for} \quad i = m+1, m+2, \ldots,
n.
\end{eqnarray}
\end{lemma}
\begin{proof} 
  Let $\eu{A}_{i} = (\phi_{i}, C)$ be ideals in $\id{O}(\rg{k})$. For
  $i > m $ these ideals are not trivial, while for $i \leq m$ they are
  equal to $\id{O}(\rg{k})$. We assume thus $i > m$ and claim that
\begin{eqnarray}
\label{qpow}
\eu{A}_{i}^{q} = (\phi_{i}).
\end{eqnarray} 
Indeed, 
\begin{eqnarray*}
\eu{A}_{i}^{q}/(\phi_{i}) =  \left(\phi_{i}^{q-1}, C \cdot \phi_{i}^{q-2},
  \ldots,  C^{q-1}, C^{q}/\phi_{i} \right).
\end{eqnarray*}
It follows from \rf{prd} that the right hand side ideal is integer and
$(\phi_{i}) \mid \eu{A}_{i}^{q}$. On the other hand, since $(\phi_{i},
\phi_{j}) = (1)$ for $i \neq j$, we have $\left(C^{q}/\phi_{i},
  \phi_{i}\right) = (1)$ and thus $\eu{A}_{i}^{q} = (\phi_{i})$, as
claimed. Furthermore, $q \nmid h(\rg{k})$ implies that the ideals
$\eu{A}_i, i > m$ must be principal. There are $\mu_{i} \in
\id{O}(\rg{k})$ such that $\eu{A}_{i} = (\mu_{i})$. It follows then
from \rf{qpow} that
\begin{eqnarray}
\label{alg}
({\mu}_{i}^{q}) & = & \eu{A}_{i}^{q} = (\phi_{i}) \quad \hbox{and} \\
\phi_{i} & = & \eta_{i} \cdot {\mu}_{i}^{q} \quad \hbox{for some}
\quad \eta_{i} \in \left(\id{O}(\rg{k})\right)^{\times}.\nonumber
\end{eqnarray}
This completes the proof of the lemma.
\end{proof}
\section{Fermat Equations and Proof of Theorem \ref{genfer}}
Suppose that $x, y, z$ is a non trivial solution of \rf{FC} and $q
\nmid h_{p} = h(\K)$. Then \rf{algq} holds by Corollary \ref{crhodef}.
This leads to a reduction of the initial Fermat - Catalan equation
\rf{FC} to a Fermat - like equation (i.e. involving only \textit{one }
prime exponent) in extension fields. It is likely that this reduction
may bring some progress in the \textit{general program } ennounced by
Darmon \cite{Da} for the solution of \rf{FC}. Indeed, in this
programatic paper, Darmon suggests that in order to solve general
cases of \rf{FC}, "\textit{one is naturally led to replace elliptic
  curves by certain 'hypergeometric Abelian varieties', so named
  because their periods are related to values of hypergeometric
  functions}". Our result shows however that there is a solid region
of the $(p,q)$ plane, in which the simple elliptic curves, albeit
defined over totally real (cyclotomic) extensions of $\Q$, can still do
the job \footnote{ I thank Jordan Ellenberg for pointing out to me
  that it is important to have Fermat equations defined over
  \textit{totally real} fields, thus opening a door to the use of
  Hilbert modular forms.  Consequently, the Fermat equations which we
  deduce here will have this property and be defined over the simplest
  totally real fields available in the context. It should be mentioned
  here, that a large variety of Fermat-like equations can be deduced
  from our results; there are reasons to believe that the ones we
  display may be the best point of departure for further
  (non-cyclotomic) investigations.}. 

The result in this direction was ennounced in Theorem \ref{genfer}, of which
we give a proof below.
\begin{proof}
  We shall treat Case I and Case II separately, using the definitions
  in \rf{cum1}, \rf{cum2}. Suppose first that $e = 0$ (Case I). Then
\[ \alpha \cdot \overline \alpha = (x+\zeta y)(x+\overline \zeta y) =
(x+y)^{2} - \mu xy = A^{2q} - \mu xy = \delta \cdot \nu^{q},
\] where $\delta = \varepsilon \cdot \overline \varepsilon$, $\nu =
\rho \cdot \overline \rho$ and $\mu = (1-\zeta)(1-\overline
\zeta)$ is the ramified prime above $p$ in $\K^{+}$. Since $p > 3$
there is at least one non trivial automorphism $\sigma \in G_{p}$; we
may apply this automorphism to the above equation and eliminate $x y$
from the resulting two identities:
\[ \left(\sigma(\mu) - \mu \right) \cdot A^{q} =
\sigma(\mu) \cdot \delta \cdot \nu^{q} -
\mu \cdot \sigma\left(\delta \cdot
  \nu^{q}\right).\] If $\sigma(\zeta+\overline \zeta) =
\zeta^{c}+\overline \zeta^{c}$, we note that
\[ \mu-\sigma(\mu) = (\zeta-\zeta^{c}) +
\overline{\zeta-\zeta^{c}} = (1-\zeta^{c-1}) \cdot (\zeta - \overline
\zeta^{c}) = \delta_{1} \mu,\] with $\delta_{1} \in
\Z[\zeta+\overline \zeta]^{\times}$. After division by $\mu$ in
the previous identity, we find there are three units $\delta_{1},
\delta_{2} = -\delta \cdot \frac{\sigma(\mu)}{\mu},
\delta_{3} = \sigma(\delta)$ such that
\[ \delta_{1} A^{2q} + \delta_{2} \nu^{q} + \delta_{3} \sigma(\nu)^{q}
= 0 .\] In this Case, \rf{fereq} holds with $a, b, c$ being all units.

Suppose now that $e = 1$, so $x+y=p^{q-1} \cdot A^{q}$ and $\alpha
\cdot (1-\zeta) = x+\zeta y$. In this case,
\[ \mu \cdot \alpha \cdot \overline \alpha = (x+\zeta
y)(x+\overline \zeta y) = (x+y)^{2} - \mu xy = A^{2q} \cdot
p^{2(q-1)} - \mu xy = \mu \cdot \delta \cdot \nu^{q}, \]
where again $\delta = \varepsilon \cdot \overline \varepsilon$ and
$\nu = \rho \cdot \overline \rho$. We eliminate, like previously, the
term in $x y$, thus obtaining:
\[ \left(\sigma(\mu) - \mu \right) \cdot A^{2q} p^{2(q-1)} =
\left(\mu \cdot \sigma(\mu)\right) \cdot \left( \delta \cdot
  \nu^{q} - \sigma\left(\delta \cdot \nu^{q} \right) \right), \] and
with the same $\delta_{1}$ as above, upon division by $\mu \cdot
\sigma(\mu)$,
\[ -\delta_{1} A^{2q} \cdot \frac{p^{2(q-1)}}{\sigma(\mu)} = \delta
\cdot \nu^{q} - \sigma\left(\delta \cdot \nu^{q} \right) .\] In this
Case, \rf{fereq} holds with $b, c$ being units, while $a = \delta_{1}
\cdot \frac{p^{2(q-1)}}{\sigma(\mu)}$, so $(a)$ is a power of the
ramified prime above $p$. This completes the proof of the Theorem.
\end{proof}

In the Second Case, one may wish a Fermat equation with \textit{all -
units} coefficients. This can be achieved at the cost of imposing $p
\geq 7$ and the fact that all three unknowns will be non -
rational. With this one has the following
\begin{proposition}
In the premises of Theorem \ref{genfer} and assuming that $e = 1$ (the
Second Case, thus), let $p \geq 7$ and $\K = \Q(\zeta_p)^+, \rg{A} =
\id{O}(\K)$. Then for any $\sigma \in \Gal(\K/\Q)$, there are three
units $\varepsilon_j \in \rg{A}$ and a $\nu \in \rg{A}$, such that the
ecuation:
\begin{eqnarray}
\label{fcase2}
    \varepsilon_1 \cdot X^q +  \varepsilon_2 Y^q + \varepsilon_3 Z^q = 0
\end{eqnarray}
has the solution $(X, Y, Z) = (\nu, \sigma(\nu), \sigma^2(\nu)) \in
\rg{A}^3$.
\end{proposition}
\begin{proof}
We start form the identity 
\[ \delta_{1} A^{2q} \cdot \frac{p^{2(q-1)}}{\sigma(\mu)} = \delta
\cdot \nu^{q} - \sigma\left(\delta \cdot \nu^{q} \right) \] derived
above for this Case, in the proof of the Theorem \ref{genfer}. We
shall need precise information about the units, and thus trace them
back in the proof. Let $\sigma$ be fixed and $\chi = \mu^{\sigma-1}$;
then $\delta_1 = \chi - 1$ is also a unit. The unit $\delta$ is fixed
by its $q$ - adic expansion, but we shall not require more detail here.

Now apply $\sigma$ to the previous identity and use the definition of $\chi$:
\begin{eqnarray*}
(\chi -1 )^{-1} \left(\delta \cdot \nu^{q} - \sigma(\delta \cdot
\nu^{q}) \right) & = & C/\sigma(\mu) \\ 
\sigma \left((\chi -1 )^{-1} \left(\delta \cdot \nu^{q} \right)-
\sigma\left(\delta \cdot \nu^{q} \right) \right) & = &
C/\sigma(\mu) \times \sigma(\chi^{-1}),
\end{eqnarray*}
and notice the crucial identity among units:
\begin{eqnarray}
\label{unitlc}
 \frac{1}{\chi - 1} - \sigma\left(\frac{\chi}{\chi-1}\right) = 
\frac{1}{\sigma(\chi^{-1}) - 1} =: \Delta \in \rg{A}^{\times} .
\end{eqnarray}
Thus, a linear combination of the previous equations yields:
\[  \varepsilon_1 \nu^q + \varepsilon_2 \sigma(\nu)^q + 
\varepsilon_3 \sigma^2(\nu)^q = 0, \]
with the units:
\[ \varepsilon_1 = \frac{\delta}{\chi-1}, \quad \varepsilon_2 = 
\frac{\sigma(\chi \cdot \delta)}{\sigma(\chi) - 1} \quad \varepsilon_3 = 
\sigma\left(\frac{\sigma(\delta)}{\chi-1}\right) . \]
This completes the proof. 
\end{proof}
\begin{remark}
One may use above two distinct Galois actions $\sigma, \tau$ rather
then just one and its square. The corresponding linear combination of
units in \rf{unitlc} remains a unit and one thus obtains a more
general Fermat equation with conjugate solutions.
\end{remark}
\subsection{The Case $p = 3$}
According to Beukers \cite{Be}, this case has been solved for $n = 4,
5, 17 \leq n \leq 10000$ by N. Bruin \cite{Br2} and A. Kraus
\cite{Kr}, respectively (note that here, composite exponents are taken
into consideration too). Since this leaves the general case open, it
may be interesting to deduce the associated Fermat equations. They are
given by the following:
\begin{proposition}
\label{fer3}
Let $q > 3$ be a prime for which the equation $x^{3}+y^{3} = z^{q}$
has non trivial coprime solutions in the integers. If $\K =
\Q(\sqrt{-3})$ is the third cyclotomic field, $\rg{E} \subset \K$ are
the Eisentstein integers and $\rho \in \rg{E}$ is a third root of
unity, then there is a $\beta \in \rg{E}$ such that one of the
following alternatives hold:
\begin{eqnarray}
\label{caseI3}
    \beta^{q} + \overline \beta^{q} & = & A^{q},\\ 
\label{caseII3}
    3(\rho-\rho^{2}) \cdot \left( \beta^{q} - \overline \beta^{q}\right) & =
    & A^{q},
\end{eqnarray}
where $A \in \Z$.
\end{proposition}
\begin{proof}
  The alternative above corresponds to the two cases of \rf{FC}. In
  the first case, $x+y=A^{q}$ and $\eu{A} = (x+\rho y, z)$ is a
  principal ideal. Since all the units of $\rg{E}$ are (sixth) roots
  of unity, and thus \nth{q} powers, it follows that there is a $\beta
  \in \rg{E}$ such that $\beta^{q} = \rho x + \overline \rho y$. Then
  $\beta^{q}+\overline \beta^{q} = (\rho+\overline \rho) (x+y) =
  -A^{q}$, which proves \rf{caseI3}. In the second case $x+y =
  A^{q}/3$ and \rf{caseII3} follows by a similar computation, with
  details left to the reader.
\end{proof}
It is also useful to know that the case $q | z$ can be ruled out in
both \rf{caseI3} and \rf{caseII3} by using a generalized form of
Kummer descent. We shall give details for $p > 3$ in a later section,
leaving this case as an open remark for the reader.
\section{Consequences of Class Field Theory}
In this section we shall deduce some consequences of class number
conditions in the cyclotomic fields of our interest. Of most value for
our investigation, these conditions give some control on local
properties of units of the \nth{pq} field and its subfields.

\subsection{Primary Numbers and Reflection}
We shall be interested in the sequel in the algebraic integers of the
field $\KL = \Q(\zeta, \xi)$ and its subfields. Let $\rg{R}$ be one of
the rings of integers $\Z[\zeta], \Z[\xi], \Z[\zeta, \xi]$ and
$\rg{K}$ its quotient-field. An element $\alpha \in \rg{R}$ is $q$ -
\textbf{primary} if it is a $q$ - adic \nth{q} power. The $q$ -
primary numbers build a subring $\rg{R}_{q} \subset \rg{R}$ and it is
immediately verified that $\rg{R}^{q} \subset \rg{R}_{q}$. If
$E(\rg{R}) = E = \rg{R}^{\times} \subset \rg{R}$ are the units of the
respective field, then we write $E_{q} = E \cap \rg{R}_{q}$. The ring
$\rg{R}_{q}$ induces the following equivalence relation:
\begin{eqnarray}
\label{eqvq}
\alpha =_{q} \beta \quad \Leftrightarrow \quad \exists \ \mu, \nu \in
\rg{R}_{q} \ : \ \mu \cdot \alpha = \nu \cdot \beta. 
\end{eqnarray}
If $\eu{q} = (q)$ or $\eu{q} = (1-\xi)$, depending on whether $\xi
\not \in \rg{R}$ or $\xi \in \rg{R}$, we let $S = \rg{R} \setminus
(\eu{q})$ and $\rg{R}' = S^{-1} \rg{R}$ be the corresponding
localization. One may extend the definition of $q$ - primary numbers
to $\rg{R}'$, thus obtaining the ring $\rg{R}'_{q} \subset \rg{R}'$.
The equivalence relation in \rf{eqvq} may then also be written as
\[ \alpha =_{q} \beta  \ \Leftrightarrow \ \alpha = \gamma \cdot \beta
\quad \hbox{ for some } \gamma \in \rg{R}'_{q}. \]

A number $\alpha \in \rg{R}'$ is called $q$ - \textbf{singular} if
there is a \textit{non - principal} ideal $\eu{B} \subset \K$ such
that $(\alpha) = \eu{B}^{q}$ as ideals. Let $\rg{R}"_{q}$ be the ring
of the $q$ - primary numbers which are also singular. The degenerate
case $\eu{B} = \rg{R}$ suggests allowing $E_{q} \subset \rg{R}''_{q}$.
By class field theory, the singular primary numbers $\alpha \in
\rg{R}''_{q}$ have the property that the extension $\Q(\zeta,
\xi)[\alpha^{1/q}]$ is unramified Abelian. There is thus, by Hilbert's
Theorem 94, an ideal of order $q$ of $\Q(\zeta, \xi)$ which
capitulates in this extension (see e.g. \cite{Wa}, Exercise 9.3).

We now define the number $h(p, q), h_{pq}^-$ ennounced in the
introduction.  The definition involves an explicit use of Leopoldt's
reflection theorem (see e.g. \cite{Lo}, \cite{Mi2}). The number $h(p,
q)$ will be defined so that the condition $\left(q, h(p, q)\right) =
1$ becomes the tightest \textit{easily computable} condition which
implies $q \nmid h_{p}$.  Let $X = \{ \chi_{0} : G_{p} \rightarrow
\overline \F_{q} \} $ be the set of Dirichlet characters of order $q$
and $\omega_{0} : G_{q} \rightarrow \F_{q}$ the Teichm\"uller
character. Then $\psi_{0} = \chi_{0}^{-1} \cdot \omega_{0}$ is a well
defined Dirichlet character of $G_{pq}$ in $\overline \F_{q}$ which
corresponds by reflection to $\chi$. If $\rg{k} = \F_{q}[\Im(X)]$ in
the obvious sense (with $\Im(\chi)$ being the image of the character
$\chi$ in $\Q(\zeta_{q-1})$), then one may chose a subfield $\rg{K}
\subset \Q(\zeta_{q-1})$ of the \nth{q-1} cyclotomic extension and an
integer ideal $\eu{Q}$ in this field, in such a way that $\rg{k} =
\id{O}(\rg{K})/\eu{Q}$. This allows to lift $\psi_{0}, \chi_{0}$ and
$\omega_{0}$ to characters $\chi, \psi, \omega$ with images in
$\rg{k}$ \cite{Mi2}. In particular, the generalized Bernoulli number
$B_{1, \psi}$ is given by \cite{Wa}:
\[ B_{1, \psi} = \frac{1}{(p-1)(q-1)} \sum_{(a,pq) = 1; \
  0 < a < p q } \ a \psi^{-1} \gamma_{a}  , \] 
where $\gamma_{a} \in G$ with $\gamma_{a}(\zeta \xi) = (\zeta
\xi)^{a}$. With this, we define
\begin{eqnarray}
\label{hpqdef}
B_{\omega} & = & \prod_{\chi_{0} \in X} \  B_{1, \chi^{-1} \omega}
\quad \hbox{ and } \nonumber \\
h(p, q) = h_{p}^{-} \cdot B_{\omega}.
\end{eqnarray}
Note that the Bernoulli numbers can be computed explicitly and in
general $B_{1,\psi} \not \in \Z$, but $B_{\omega} \in \Z$, since it is
the norm of an algebraic integer in $\Q(\zeta_{q-1})$. It is also true
\cite{Mi2}, that $q \mid h_{p}^{+}$ implies $q \nmid B_{\omega}$.

We shall see that certain ideals of $\KL$ occurring in the subsequent
proofs have order dividing $q$; they are thus principal is $q \nmid
h_{pq}^{-}$ (by reflection then $q \nmid h_{pq}$!); if the ideals
belong to $\K$, then they are already principal if $q \nmid h(p, q)$.
The following Proposition reflects these and further useful
consequences of the above class number conditions.

\begin{proposition}
\label{lemmaA}
Let $E, E_{q} \subset \KL$ be the units, respectively the $q$ -
primary units of the \nth{pq} cyclotomic extension and $E', E'_{q}
\subset \K$ be the respective sets in the \nth{p} cyclotomic
extension. If $q \nmid h_{pq}^-$ then $E_{q}=E^{q}$ and in particular,
if $\varepsilon =_{q} 1$ is a unit, then it is a \nth{q} power.
Likewise, if $q \nmid h(p, q)$, then $E'_{q} = {E'}^{q}$ and all $q$ -
primary units in $E'$ are \nth{q} powers.  Furthermore, $q \nmid
h_{pq}$ in the first case and $q \nmid h_{p}$ in the second.
\label{qprim}
\end{proposition}
\begin{proof}
  We start with the implication $q \nmid h_{pq}^{-} \Rightarrow q \nmid
  h_{pq}$; this follows directly by reflection in $\KL$. Let $q \nmid
  h_{pq}^-$ and $\varepsilon \in E_{q} \setminus E^{q}$ be a $q$ -
  primary unit, which is not a \nth{q} power. Then $\rg{k} =
  \LK(\varepsilon^{1/q})$ is an Abelian unramified extension (see e.g
  \cite{Wa}, lemma 9.1, 9.2) and $[\rg{k}:\KL] = q$; by Hilbert's
  Theorem 94, there is an ideal of order $q$ from $\KL$ which
  capitulates in $\rg{k}$, in contradiction with $q \nmid h_{pq}^{-}$.
  
  We now consider the case $q \nmid h(p, q)$. The crucial remark here
  is that even Dirichlet characters of $\K$ correspond by reflection
  to the characters indexing Bernoulli numbers which divide
  $B_{\omega}$ above. The claims for this case follow then in analogy
  to the ones for $q \nmid h_{pq}^{-}$.
\end{proof}

\subsection{Units}
We start the analysis of local properties of units in $\id{O}(\KL)$
under certain eventual restriction on the class number, by several
simple basic Lemmata.
\begin{lemma}
\label{lru}
Let $\delta = < -\zeta \cdot \xi >$ be a root of unity of $\KL$. If
$\delta \equiv 1 \mod q$, then $\delta = 1$. 

Furthermore, if $\varepsilon \in \Z[\zeta, \xi]^{\times}$ is a unit
such that $\varepsilon = a + b(\zeta) q + O(q\lambda)$ and $a \in \Z$
and $b(\zeta) \in \Z[\zeta]$, then $b(\zeta) = \overline b(\zeta)$.
\end{lemma} 
\begin{proof}
Let $\delta^2 = \zeta^a \cdot \xi^b$, with $0 \leq a < p, \ 0 \leq b
< q$ - squaring cancels the sign. If $a \cdot b \neq 0$, then
$\delta^2$ is a primitive \nth{pq} root of unity and thus
\[ P(G) = \prod_{\psi \in G} \left(1-\psi(\delta) \right) = 
\Phi_{pq}(1) = 1.\] But since $\delta^2 \equiv 1 \mod q$ we should
have $P(G) \equiv 0 \mod q^{\varphi(pq)}$, so $a \cdot b \neq 0$ is
impossible. In the cases $a = 0, b \neq 0$ and $a \neq 0, b= 0$, the
root $\delta^2$ is primitive of order $q$, resp. $p$ and the value of
$P(G)$ is $q^{p-1}$ and $p^{q-1}$, respectively. In both cases $P(G)
\not \equiv 0 \mod q^{\varphi(pq)}$, so we must have plainly $\delta^2 = 1$
and since $\delta \equiv 1 \mod q$, also $\delta = 1$.

If $\varepsilon$ is like in the claim of the lemma, then $\delta =
\varepsilon/\overline \varepsilon \equiv 1 \mod q$ is a root of unity,
and thus $\varepsilon = \overline \varepsilon$. The claim on
$b(\zeta)$ follows.
\end{proof}

\begin{lemma}
\label{unitp}
Let the $p, q$ be odd primes with $q \nmid h_{pq}^{-}. p \nmid (q-1)$
and suppose that $\varepsilon \in \id{O}(\KL)^{\times}$ is a unit such
that $\varepsilon =_{q} a$, with $a \in \Z$. Then $\varepsilon$ is a
\nth{q} power.
\end{lemma}
\begin{proof}
  Let $\sigma \in G = \Gal(\Q(\zeta, \xi)/\Q)$. Then $\delta_{\sigma}
  = \varepsilon^{1-\sigma} =_{q} 1$ and by Proposition \ref{lemmaA} it
  is a \nth{q} power. Since this holds for all $\sigma$,
\[ \varepsilon^{(p-1)(q-1)} = \prod_{\sigma \in G} \delta_{\sigma} \in
  \KL^{q} . \] But $\left(q, (p-1)(q-1)\right) = 1$ by hypothesis and
  consequently $\varepsilon$ in a \nth{q} power, as claimed.
\end{proof}

Finally, we have:
\begin{lemma}
\label{lemmaC}
\label{lamunit}
Let $p, q$ be primes and $\varepsilon \in \Z[\zeta_{p}]^{*}$ be a unit
such that $\varepsilon =_{q} c(1-\zeta)$, with $c \in \Q$. If $p \not
\equiv 1 \mod q$ and $(p-1) m \equiv 1 \mod q$, then $p =_{q} 1$ and
\begin{eqnarray}
\label{uc}
 \varepsilon =_{q} \left(\frac{1-\zeta}{p^{m}}\right) =_{q}\left(
  \frac{(1-\zeta)^{p-1}}{p}\right)^{m} = \gamma \in \Z[\zeta]^{\times} .
\end{eqnarray}
\end{lemma}
\begin{proof}
  Let $\sigma \in G = \Gal(\Q(\zeta_{p})/\Q)$ be a generator and
  $\Omega = \sum_{i=0}^{p-3} (p-2-i) \sigma^{i} \in \Z[G]$, so that
  $(\sigma-1) \Omega + (p -1)= \Norm_{\Q(\zeta)/\Q}$. By hypothesis,
  $\varepsilon^{\sigma-1} =_{q} \eta = (1-\zeta)^{\sigma-1}$, so
\[ 1 = \Norm(\varepsilon) = \varepsilon^{p-1+\Omega(\sigma-1)} =_{q}
\varepsilon^{p-1} \cdot \eta^{\Omega}, \]
and thus $\varepsilon^{p-1} =_{q} \eta^{-\Omega}$ and $\varepsilon =_{q}
\eta^{- m \Omega}$. Note also that $\eta^{\Omega} =
(1-\zeta)^{(\sigma-1)\Omega} = (1-\zeta)^{\mathbf{N}-p+1} =
p/(1-\zeta)^{p-1}$, a simple expression for this unit. 
\end{proof}

The deeper results on local properties of units in $\KL$ (and
subfields), given class number restraints, are summarized in the
following Proposition. The proof of the proposition is quite lengthy
and involves, along with the previous class field results, some
interesting properties of cyclotomic units in fields of composite
order. The result is interesting in itself, but it shall be used only
for improving an estimate in Theorem \ref{main2} for the special case
when $q \equiv 1 \mod p$. Given the lack of generality of its
application, the reader who is more interested in an overview of the
main ideas and proofs, can thus skip to the next section.
\begin{proposition}
\label{un1q}
Let $p, q$ be odd primes with $q \nmid h_{pq}^-$, $p \not \equiv 1 \mod
q$ and $\KL = \Q[\zeta, \xi]$ be the \nth{pq} cyclotomic extension of
$\Q$. If $\varepsilon \in \KL$ is a unit for which there is a $a \in
\Z[\zeta]$ such that $\varepsilon \equiv a \mod q$, then there is a
unit $\delta \in \Z[\zeta]$ such that $\varepsilon =_q \delta$.
\end{proposition}
We first prove some Lemmata:
\begin{lemma}
\label{linindep}
Let $p, q$ be odd primes with $p \not \equiv 1 \mod q$ and $q \nmid
h_p^-$; let $\zeta \in \C$ be a primitive \nth{p} root of unity. Then
$\gamma_c = 1/(1-\zeta^c), c=1, 2, \ldots, p-1$ form a basis for the
Galois ring $\Z[\zeta]/\left(q \cdot \Z[\zeta]\right)$.
\end{lemma}
\begin{proof}
  We must show that $\gamma_c$ are linear independent. For this, we
  shall show that the discriminant $\Delta$ of the $\Z$ - module $M =
  [\gamma_1, \gamma_2, \ldots, \gamma_{p-1}]$ is coprime with $q$
  under the given conditions. Let $\sigma$ be a generator of
  $G_{p}=\Gal(\Q(\zeta)/\Q)$ and the matrix $\mathbf{A} =
  \left(\sigma^{i+j} \left(\frac{1}
      {1-\zeta}\right)\right)_{i,j=1}^{p-1}$; then $\Delta =
  \det^2(\mathbf{A})$. The matrix $\mathbf{A}$ is a circulant matrix;
  if $\omega \in \C$ is a primitive \nth{p-1} root of unity, then the
  vectors $\vec{f}_k = (\omega^{jk})_{j=0}^{p-2}$, for $k = 0, 1,
  \ldots, p-2$ are eigenvectors of $\mathbf{A}$. In the base spanned
  by these vectors, $\mathbf{A}$ is diagonal. The base transform
  matrix having $\vec{f}_k$ as columns is a Vandermonde matrix with
  discriminant $D = \prod (\omega^i - \omega^j)$, which is a power of
  $p-1$ and thus a unit modulo $q$, since $p \not \equiv 1 \mod q$.
  The base transform is thus regular modulo $q$ and $\Delta$ is
  invertible modulo $q$ iff its conjugate matrix in the eigenvector
  base is so.
  
  We now compute the determinant of the conjugate diagonal matrix
  $\rg{D}$ of $\rg{A}$. If $\chi : \ZM{p}^{\times} \rightarrow <
  \omega >$ is a character with $\chi(\widehat{\sigma}) = \omega$ -
  where $\sigma(\zeta) = \zeta^{\widehat{\sigma}}$ - then one verifies
  that the representation of $\mathbf{A}$ in the eigenvector base is
\[ \mathbf{A} \sim \diag\left(\tau'(\chi^j)\right)_{j=0}^{p-2}, \]
where, for $\psi \in < \chi > $, we defined the Lagrange resolvents
\[ \tau'(\psi) = \sum_{x \in \left(\Z/p\Z\right)^{\times}} 
\frac{\psi(x)}{1-\zeta^x} .\] Let $\tau(\psi) = \sum_{x=1}^{p-1}
\psi(x) \zeta^x$ be a regular Gauss sum; a known formula (see e.g.
\cite{La}) implies $\psi^{-1}(i) \tau(\psi) = \sum_{x=1}^{p-1} \psi(x)
\zeta^{ix}$. By an easy calculation, $-p/(1-\zeta^x) =
\sum_{i=1}^{p-1} i \cdot \zeta^{ix}$. Finally, assembling these
formulae we find:
\begin{eqnarray*}
 -p \cdot \tau'(\psi) & = & \sum_{x \in \left(\Z/p\Z\right)^{\times}}
\psi(x) \cdot \sum_{i=1}^{p-1} i \cdot \zeta^{xi} \\ & = &
\sum_{i=1}^{p-1} i \cdot \left(\sum_{x \in
\left(\Z/p\Z\right)^{\times}} \psi(x) \cdot \zeta^{xi}\right) =
\tau(\psi) \cdot \sum_{i=1}^{p-1} i \cdot \psi^{-1}(i).
\end{eqnarray*}
But $B_{1,\psi^{-1}} = \frac{1}{p} \cdot \sum_{i = 1}^{p-1} i \cdot
\psi^{-1}(i)$ is a generalized Bernoulli number and
$\Norm(B_{1,\psi^{-1}}) | h_p^-$. But since $q \nmid h_{p}^{-}$, this
is a unit modulo $q$. Furthermore, $\tau(\psi) \cdot \overline
\tau(\psi) = p$, so $\tau(\psi)$ is a unit too. Finally $\tau'(\psi) =
- \tau(\psi) \cdot B_{1,\psi^{-1}}$ is a unit modulo $q$ for all $\psi
\in < \chi >$.

But 
\[\det(\mathbf{A})= \det \left(\diag
\left(\tau'(\chi^j)\right)_{j=0}^{p-2}\right) = \prod_{j=0}^{p-2} \ 
\tau'(\chi^j) . \] Since all the factors have been shown to be units
modulo $q$, it follows that $(\Delta, q) =
\left(\det^2(\mathbf{A}),q\right) = 1$, which completes, the proof.
\end{proof}

We study next the cyclotomic units of $\KL$. Let $\delta = (1 - \zeta
\xi)$ and $C_{1} = < - \ \zeta \xi > \cdot \left(\Z[G]^{+} \delta
\right) \subset \Z[\zeta, \xi]^{\times}$ be the $Z[G]$ module
generated by the unit $\delta$ together with the roots of unity. If $q
\not \equiv 1 \mod p$ we let $C_{2} = \{1\}$; otherwise, let $C_{2} =
\Z[G_{p}] \cdot (1-\zeta)^{(1+\jmath)(\sigma-1)}$ be the
$\Z[G_{p}]^{+}$ module of units of $\K^{+}$ generated by the
cyclotomic unit $\eta = |(1-\zeta)^{\sigma-1}|^{2}$. Thus $C_{2}$ is
in this case \textit{essentially} equal to the cyclotomic units of
$\K^{+}$ \cite{Wa}, Chapter VIII; it has in fact an index $2$ in this
group, which is of no importance in our context, since we are focusing
on $q$ - parts of unit groups.

Note that for $q \equiv 1 \mod p$, the norm
$\Norm_{\KL^{+}/\K^{+}}(\delta) = \frac{1-\zeta^{q}}{1-\zeta} =
1$ and we always have
\begin{eqnarray}
\label{disjoint}
C_{1} \cap C_{2} = \{1\}.
\end{eqnarray}
For $q \not \equiv 1 \mod p$, the statement is trivial. Suppose now
that $p | (q-1)$ and let $\varepsilon = \gamma_{1} \cdot \gamma_{2}
\in C_{1} \cap C_{2}$, with $\gamma_{i} \in C_{i}$. Since $\varepsilon
\in C_{1}$, $\Norm_{\KL^{+}/\K^{+}}(\varepsilon) = 1 =
\gamma_{2}^{q-1}$ and since $\gamma_{2}$ is real, it follows that
$\gamma_{2} = \pm 1$. But $\varepsilon \in C_{2}$ implies, by taking
norms again, that $(\varepsilon/\gamma_{2})^{q-1} = \gamma_{1}^{q-1} =
1$ and eventually $\varepsilon = \gamma_{1} \cdot \gamma_{2} = 1$, as
claimed.
\begin{lemma}
\label{q1p}
Let $p, q$ be odd primes with $q \nmid h_{pq}^-$ and $p \not \equiv 1
\mod q$. If $C' = C_{1}$ for $q \not \equiv 1 \mod p$ and $C' = C_{1}
\cdot C_{2}$ otherwise, then $C'$ has finite index in $E$, the group
of units of $\KL$ and $q \nmid \kappa = [E : C']$. In particular,
\begin{eqnarray}
\label{EC}
E = C' \cdot E^{q}.
\end{eqnarray}
\end{lemma}
\begin{proof}
  In the case $q \not \equiv 1 \mod p$, there are no multiplicative
  dependencies in $C_{1}$ and the claims are a direct consequence of
  \cite{Wa} Corollary 8.8 (note that both $C'$ and $E$ contain the
  same torsion, the roots of unity of $\KL$). Indeed, since $p \not
  \equiv 1 \mod q$ and $q \not \equiv 1 \mod p$, the Euler factor in
  this corollary is not vanishing and also coprime to $q$. Thus
  $E/E^{q}$ and $C'/{C'}^{q}$ have the same rank and annihilators and
  the subsequent claims follow from this observation, since $C_{2}$ is
  trivial in this case.
  
  We now consider the case $q \equiv 1 \mod p$, for which we apply the
  Theorem 8.3 in \cite{Wa}. Note that the Ramachandra units are, up to
  roots of unity and an index $4$, exactly $C' = C_{1} \cdot C_{2}$ in
  this case. By the Theorem of Ramachandra, it follows that $q \nmid
  \kappa = [ E : C']$, which also implies \rf{EC}.
\end{proof}

Next we investigate the structure of the cyclotomic units as group
ring modules. For this we note that the map $\iota : \Z[G] \rightarrow
\id{Z} = \Z[X, Y]/\left(X^{p-1} - 1, Y^{q-1}-1\right)$ given by
$\sigma \mapsto X$ and $\tau \mapsto Y$ - where $\sigma, \tau$ are
generators of $G_{p}, G_{q}$, as usual - is an isomorphism of rings.
We shall consider next various restrictions of this map to subrings
and quotient rings of $\Z[G]$, without changing the notation.

The image of $\Z[G]^{+}$ under this map is $\id{Z}^{+} =
\id{Z}/(X^{(p-1)/2} - Y^{(q-1)/2})$, since the partial conjugations
$\jmath_{p} = \jmath_{q}$ in the real subfield.

We are interested in the $q$ - parts $W' = C'/{C'}^{q}$ and the
components $W_{i} = C_{i}/C_{i}^{q}$, for $i = 1, 2$. These are
obviously $\F_{q}[G]^{+}$ - modules and as a consequence of
\rf{disjoint} we also have 
\[ W' = W_{1} \oplus W_{2}. \]
Since the $\Norm_{\KL^{+}/\Q}$ annihilates the units, they are also
$\rg{R} = \F_{q}[G]^{+}/ \left(\Norm_{\KL^{+}/\Q}\right)$ - modules.
We have
\begin{eqnarray*} 
\rg{R} & \cong & \id{R} =  \id{Z}^{+}/\left(q, \frac{X^{p-1} \cdot
    Y^{q-1} - 1} {X    \cdot Y - 1}\right) \\
& = & \F_{q}[X, Y]/\left(X^{p-1} - 1, Y^{q-1}-1,
  X^{(p-1)/2} - Y^{(q-1)/2}, \frac{X^{p-1} \cdot Y^{q-1}-1}{X
    \cdot Y - 1} \right),
\end{eqnarray*}
under the isomorphism $\iota$. The above isomorphism illustrates that
$\rg{R}$ is a semi-simple module, which is not cyclic.

Suppose now that $q \not \equiv 1 \mod p$ so $C' = C_{1}$, a cyclic
$\rg{R}$ - module. By comparing ranks in \rf{EC}, it follows in fact
that $C_{1} = \rg{R} \cdot \delta$ in this case. If $q \equiv 1 \mod
p$ then $\Norm_{\KL^{+}/\K^{+}}(\delta) = 1$ yields some
multiplicative dependencies in $C_{1}$. If $\rg{R}_{1} =
\rg{R}/\left(\Norm_{\KL^{+}/\K^{+}}\right) \cong
\iota\left(\rg{R}\right)/(X^{(q-1)/2}-1)$, then one verifies that
$W_{1} = \rg{R}_{1} \cdot \delta$ in this case; in order to keep a
uniform notation, we shall also write $\rg{R}_{1} = \rg{R}$, if $q
\not \equiv 1 \mod p$, so that $W_{1} = \rg{R}_{1} \cdot \delta$ in
both cases.

As to $W_{2}$, by \cite{Wa}, Theorem 8.11, one simply has $W_{2} =
\F_{q}[G_{p}]^{+} \cdot \eta$. Note that the ranks of $\rg{R}_{1}$ and
$\F_{q}[G_{p}]^{+}$ add up to $(p-1)(q-1)/2 = \rk{W'}$. We now apply
the gained structure for analyzing some particular cyclotomic units.
\begin{lemma}
\label{ltau}
  The notations being like above, let $\delta_{1} \in C_{1}$ be a unit
  which verifies $\delta_{1}^{\tau-\hat{\tau}} \in C_{1}^{q}$. Then
  $\delta_{1}^{\sigma \tau -1} = \delta^{\Theta}$ for some $\Theta \in
  \rg{R}_{1}$ such that
\begin{eqnarray}
\label{annih1}
\Theta = \varepsilon_{1} \cdot \theta, \quad \hbox{with} \quad  \theta \in  
\F_{q}[G_{p}]^{+},
\end{eqnarray}
and $\varepsilon_{1} = \frac{1}{q-1} \sum_{b=1}^{q-1} b \cdot
\tau_{b}^{-1} \in G_{q}$ is the first orthogonal idempotent of
$\F_{q}[G_{q}]$. Furthermore,
\begin{eqnarray}
\label{taum1}
\delta_{1} \in C_{1} \quad \hbox{and} \quad \delta_{1}^{\tau-1} \in
C_{1}^{q} \quad \Rightarrow \quad \delta_{1} \cdot \Norm(\delta_{1}) \in C_{1}.
\end{eqnarray}
\end{lemma}
\begin{proof}
  We let $\widetilde{\delta}_{1} = \delta_{1} \mod C_{1}^{q}$ be the
  image in $W_{1}$. Then the hypothesis on $\delta_{1}$ translates to
  $\widetilde{\delta}_{1}^{\tau-\hat{\tau}} = 1$. If $\delta_{1} =
  \delta^{\Theta_{0}}$ for some $\Theta_{0} \in \Z[G]^{+}$, then
  $\Theta(\tau-\hat{\tau})$ lays thus in the kernel of the map
  $\Z[G]^{+} \rightarrow \rg{R}_{1}$.
  
  We shall have, like usual, to distinguish whether $q \equiv 1 \mod
  p$ or not.  In the latter, simple case, we know that $W_{1} = \rg{R}
  \cdot \delta$ and the previous remark on $\Theta$ implies that 
\[ \iota(\Theta_{0}) \cdot (\tau-\hat{\tau}) \equiv 0 \mod \left(q,
  \frac{X^{p-1} \cdot Y^{q-1}-1}{X \cdot Y - 1}  \right).\] In the
second case, we have
\[ \iota(\Theta_{0}) \cdot (\tau-\hat{\tau}) \equiv 0 \mod \left(q,
  \frac{Y^{q-1}-1}{Y - 1} \right).\] We let $\Theta = \Theta_{0} \cdot
(\sigma \tau -1) $.  The second generators of the ideals in the
kernels of the last two congruences are images of norms and they are
annihilated in $\id{Z}^{+}$ by $\iota(\sigma \tau-1)$. It follows that
$\Theta \cdot (\tau-\hat{\tau}) \equiv 0 \mod q$.  Let $\theta_i \in
\Z[G_p], i = 0, 1 , 2, \ldots, q-2$ be such that
\begin{eqnarray*}
 \Theta & = & \sum_{n=1}^{q-1} \tau^n \cdot \theta_{n-1} \quad \hbox{and}
\\ \Theta \cdot (\tau-\widehat{\tau}) & = & \sum_{n=1}^{q-1}
\left(\theta_{n-2} - \widehat{\tau} \theta_{n-1} \right) \tau^{n} \equiv 0
\mod q,
\end{eqnarray*}
where the indices in the last sum are taken modulo $q-1$. Since
${\tau^n}$ are independent over $\Z[G_p]$, the sum vanishes modulo
$q$ if all of the coefficients do. Thus, inductively,
\[  \theta_n \equiv \left(\widehat{\tau}\right)^{-n} \cdot \theta_0, 
\quad n = 1, 2, \ldots, q-2. \]
But then 
\[\Theta \equiv -\theta_0 \cdot \varepsilon_1 \mod q\Z[G^+], \]
with $\varepsilon_1 \equiv -\sum_{n=1}^{q-1} (\tau/\widehat{\tau})^n =
-\sum_{a=1}^{q-1} a \tau_a^{-1} \mod q$ being the first orthogonal
idempotent of $\ZM{q}[G_q]$. 

We now prove \rf{taum1}. For this we note the following decomposition
in $\Z[G_{q}]$: $\Norm = \Norm_{\K'/\Q} = \Norm_{\KL/\K} = (\tau-1)
\cdot \Omega + (q-1)$, for some $\Omega \in \Z[G_{q}]$; the
verification is a simple computation and is left to the reader. But
then, given $\delta_{1}$ in \rf{taum1}, we have:
\[ \Norm(\delta_{1}) = \delta_{1}^{(\tau-1) \Omega + q}/\delta_{1}
\quad \hbox{and} \quad \delta_{1} \cdot \Norm(\delta_{1})
\in C_{1}^{q}. \] This completes the proof.
\end{proof}

The main result towards the proof of the Proposition is the following:
\begin{lemma}
\label{unitq1}
Let $p, q$ be odd primes with $q \nmid h_{pq}^-$, $p \not \equiv 1
\mod q$, $\KL= \Q(\zeta, \xi)$ be the \nth{pq} cyclotomic extension and
$G = \Gal(\KL/\Q) = G_p \times G_q$ with $G_p = \Gal(\Q(\zeta)/\Q), G_q
= \Gal(\Q(\xi)/\Q)$. If $\varepsilon \equiv 1 \mod q \lambda$ is a
unit of $\KL$, then $\varepsilon$ is a \nth{q} power. In particular,
if $\varepsilon = 1 + a q \lambda + O(q \lambda^{2})$, with $a \in
\Z[\zeta]$, then there is a $\beta \in \Z[\zeta]$ such that
\begin{eqnarray}
\label{addh90}
a \equiv \sigma_{q}(\beta) - \beta \mod q.
\end{eqnarray}
\end{lemma} 
\begin{proof}
  Let $E \subset \KL$ be the real units, $C' \subset E$ the cyclotomic
  units defined above and $\varepsilon \equiv 1 + a q \lambda \mod q
  \lambda^2$, for $a \in \Z[\zeta]$. If $\kappa = [E : C'] \in \N$,
  then $\varepsilon^{\kappa} \in C'$ is a unit with the same type of
  $\lambda$ - expansion as $\varepsilon$, since $(\kappa, q) = 1$. We
  may thus assume, for simplicity, that $\varepsilon \in C'$ to start
  with and thus
\[ \varepsilon = \delta_{1} \cdot \delta_{2}, \quad \hbox{with} \quad
\delta_{i} \in C_{i}. \] 

Note that $\varepsilon^{\tau-1} = \delta' = \delta_{1}^{\tau-1} \equiv
1 \mod q \lambda$. Since $\tau(\lambda) \equiv \widehat{\tau} \cdot
\lambda \mod \lambda^2$, we have
\[\psi = {\delta'}^{\tau-\hat{\tau}} = \varepsilon^{(\tau-1)(\tau- 
  \widehat{\tau})} \equiv \frac{1 + a' \cdot q \cdot \tau(\lambda)} {1
  + a' \cdot q \cdot \widehat{\tau} \cdot \lambda} \cdot \lambda
\equiv 1 \mod q \lambda^2, \] where $a' = \widehat{\tau}-1 \in
\Z^{\times}$; thus $\psi=_q 1$ and so $\psi \in {C'}^q \cap C_{1} =
{C_{1}}^{q}$ by Proposition \ref{lemmaA} and \rf{EC}.

We are thus in the context of Lemma \ref{ltau}, which implies that 
\[ \delta' = \delta^{\ - \varepsilon_{1} \cdot \theta} \equiv 1 \mod q
\lambda.\] We now estimate the unit $\delta^{\varepsilon_{1}}$ up to
$\lambda^2$ and compare the result with the above.
\begin{eqnarray*}
\delta & = & 1-\zeta \xi = (1-\zeta) + \zeta(1-\xi) = (1-\zeta) \cdot
\left(1 - \frac{1-\xi}{1-\overline \zeta}\right), \quad \hbox{so} \\
\delta^{-\varepsilon_1} & \equiv & (1-\zeta)^{-\varepsilon_1} \cdot
\left(1 -\sum_{a=1}^{q-1} a \cdot \tau_a^{-1}\left(\frac{1-\xi}{1-\overline
\zeta}\right) + O(\lambda^2)\right) \\ 
& \equiv &
(1-\zeta)^{-q(q-1)/2} \cdot \left(1+\frac{1-\xi}{1-\overline \zeta} +
O(\lambda^2)\right) \ \ \mod q \Z[\zeta, \xi].  
\end{eqnarray*}

If $\theta = \sum_{c=1}^{p-1} n_c \sigma_c$ and $A =
(1-\zeta)^{-q(q-1)/2} \in \Z[\zeta]^{q}$, then
\begin{eqnarray*}
\delta^{-\varepsilon_1} & = & A \left(1 +\frac{1-\xi}{1-\overline \zeta} +
O(\lambda^2)\right) \quad \hbox{and} \\ 
\delta' = \delta^{\Theta} & = & \delta^{-\varepsilon_1 \theta} = 
A^{\theta} \left(1+\frac{1-\xi}{1-\overline \zeta} +
O(\lambda^2)\right)^{\theta} \\ & = & A^{\theta} \left(1 + (1-\xi) \cdot  
\sum_{c=1}^{p-1} \frac{n_c}{1-\zeta^{-c}} + O(\lambda)^2\right).
\end{eqnarray*}
Lemma \ref{linindep} implies that the sum in the last equation only
vanishes modulo $q$ if all the coefficients $n_c$ vanish, so $\theta
\equiv 0 \mod q$ and also $\Theta \equiv 0 \mod q$, to start with.
But then $\delta_{1}^{\tau-1} = \delta' = \delta^{\Theta}$ is a
\nth{q} power. If $q \equiv 1 \mod p$, then \rf{taum1} implies that
$\delta_{1}$ is a \nth{q} power, and since we have already shown that
$\delta_{2}$ is a \nth{q} power, we have $\varepsilon = \delta_{1}
\delta_{2} \in E^{q}$. Oddly, the case $q \not \equiv 1 \mod p$
requires now more attention - this is not an intrinsic problem, but
rather a consequence of the build up of the auxiliary Lemmata, where
the load was taken away from the second case. For the case $q \not
\equiv 1 \mod p$ we have thus, again using \rf{taum1}, that
$\varepsilon = \psi \cdot \gamma^{q}$, with $\psi =
\Norm(\delta_{1})^{-1} \in \Z[\zeta]$ and $\gamma \in C_{1} = C'$. If
$\gamma = a + b \lambda$, with $a \in \Z[\zeta]$ and $b \in \Z[\zeta,
\xi]$, then the definition of $\varepsilon$ implies that
\[ \psi \equiv a^{-q} \mod q \lambda ,\]
and since $x \equiv 0 \mod \lambda$ implies $x \equiv 0 \mod q$ for $x
\in \Z[\zeta]$, it follows that $\psi =_{q} 1$ and by Proposition
\ref{lemmaA} it follows that $\psi$ is a \nth{q} power. 

We still have to prove \rf{addh90}. Let $\gamma^{q} = \varepsilon = 1
+ a q \lambda + O(q\lambda^{2})$. Then $\gamma = 1 + b \lambda +
O(\lambda^{2})$, with $b \in \Z[\zeta]$ and raising to the \nth{q}
power we find that 
\[ a \equiv (b^{q} \lambda^{q-1}/q + b) \mod \lambda .\]
But $\frac{\lambda^{q-1}}{q} = \prod_{i=1}^{q-1}
\frac{1-\xi}{1-\xi^{i}} \equiv \prod_{i=1}^{q-1} (1/i) \equiv -1 \mod
\lambda$, where the last congruence is derived from Wilson's theorem.
Thus $a \equiv b^{q} - b \mod \lambda$ and since $a, b \in \Z[\zeta]$,
it follows that $a \equiv b^{q} - b \mod q$. If $\eu{Q}$ is a prime of
$\Z[\zeta]$ over $q$, then it is fixed by $\sigma_{q}$ and $b^{q}
\equiv \sigma_{q}(b) \mod \eu{Q}$, so the previous equivalence implies
$a \equiv \sigma_{q}(b) - b \mod \eu{Q}$. This holds for all primes
$\eu{Q} | (q)$ uniformly, and it follows that $a \equiv \sigma_{q}(b)
- b \mod q$. This completes the proof.
\end{proof}

We can now prove Proposition \ref{un1q}:
\begin{proof}
  By the hypothesis of the Proposition, one can write $\varepsilon = a
  + b \cdot q + q \lambda \nu$, with $a, b \in \Z[\zeta]$ and $\nu \in
  \KL$; then
\[ \varepsilon = C \cdot (1 + q \lambda \nu'),  \quad \hbox{with } \quad
C = a + b \cdot q \in \Z[\zeta], \ \nu' = \nu/C. \] We let
$\delta^{-1} = \Norm_{\KL/\K}( \varepsilon) = C^{q-1} \cdot
\Norm_{\KL/\K} (1+ q\lambda \nu') =_q C^{-1}$ and thus $\delta =_q C
\equiv \varepsilon \mod q \lambda$. Obviously, $\delta \in
\Z[\zeta]^{\times}$ and $\varepsilon/\delta = 1+q \lambda \nu''$ is a
unit verifying the hypothesis of Lemma \ref{unitq1}. The first claim
follows by applying the Lemma to $\varepsilon/\delta$ (note that the
claim is trivial if $\varepsilon \in \Z[\zeta]$).

\end{proof}

\section{An Improved Case Distinctions}
In this section we derive some easy consequences from the conditions
deduced in the previous one. Finally, the methods developed in this
section will be sharpened in the next one, thus leading to a proof of
Theorem 1. In this Theorem, the two Cases discussed above, and which
depend on congruences modulo $p$, analogous to the Abel-Barlow Cases
in the classical Fermat equation, split into three additional cases
each, and these additional cases rely upon congruences modulo $q$.
\begin{lemma}
\label{tdev}
  Let $p, q$ be odd primes and $x, y$ coprime integers with $x \cdot y
  \not \equiv 0 \mod q$ and such that there is a $\beta \in \Q(\zeta)$
  with
\begin{eqnarray}
\label{premq}
 \frac{x+\zeta^{q} \cdot y}{x+\overline \zeta^{q} \cdot y} = \pm
\left(\frac{\beta}{\overline \beta}\right)^{q} .
\end{eqnarray}
Then
\begin{eqnarray}
\label{ponder} -(\zeta^{q}-\overline \zeta^{q}) \varphi(t)  \equiv 
 \sum_{k=1}^{q-1}  \ \frac{t^{k}-t^{2-k}}{k} \cdot 
(\zeta^{k}-\overline \zeta^{k}) \mod q. 
\end{eqnarray}
\end{lemma}
\begin{proof}
  A development of \rf{premq} up to the second power of $q$ yields:
\[ \frac{x+\zeta^{q} \cdot y}{x+\overline \zeta^{q} \cdot y} \equiv 
\left(\frac{x+\zeta \cdot y}{x+\overline \zeta \cdot y}\right)^{q}
\mod q \Z[\zeta] . \]
Combining with \rf{premq} we find 
\[ \pm \frac{\beta}{\overline \beta} = \frac{x+\zeta \cdot
  y}{x+\overline \zeta \cdot y} + q \cdot \mu, \]
with $\mu \in \Q(\zeta)$ being a $q$ - adic integer. Raising to the power
$q$, it follows that in fact 
\begin{eqnarray}
\label{q2}
 \frac{x+\zeta^{q} \cdot y}{x+\overline \zeta^{q} \cdot y} \equiv 
\pm \left(\frac{x+\zeta \cdot y}{x+\overline \zeta \cdot y}\right)^{q}
\mod q^{2} \Z[\zeta] . 
\end{eqnarray}

We write $\varphi(a) = \frac{a^{q}-a}{q} \mod q$, for $(a,q) = 1$ and
let $t \equiv -y/x \mod q^{2}$, so $-(y/x)^{q} \equiv t + q \varphi(t) \mod
q^{2}$. Now 
\begin{eqnarray*} 
(x+\zeta \cdot y)^{q} & \equiv & x^{q} \cdot (1 - t\cdot \zeta)^{q}
 \equiv  (x+q \varphi(x)) \cdot
(1-t\zeta)^{q} \\ & \equiv & (x+q \varphi(x)) \cdot \left(1-t\zeta^{q} +  q
  f(\zeta)\right) \mod q^{2} \quad \hbox{where} \\
f(\zeta) & = & -\zeta^{q} \cdot \varphi(t) + \sum_{k=1}^{q-1} 
\binom{q}{k} (-t\zeta)^{k} \equiv - \left(\zeta^{q} \cdot \varphi(t) +   
\sum_{k=1}^{q-1} \frac{ t^{k} \zeta^{k}}{k} \right) \mod q.
\end{eqnarray*}

Writing $x+\zeta^{q} y = x(1-t\zeta^{q}) = x \cdot \alpha$ and
eliminating denominators in \rf{q2} we find that
\begin{eqnarray*}
\alpha \cdot (x+\varphi(x)) \left(\overline \alpha + q \cdot f(\overline
  \zeta)\right) & \equiv & \overline \alpha \cdot (x+\varphi(x))
  \cdot  \left(\alpha + q \cdot
  f(\zeta) \right) \mod q^{2} \quad \hbox{and} \\
\alpha \cdot  f(\overline  \zeta) & \equiv & \overline \alpha \cdot
  f(\zeta) \mod q .
\end{eqnarray*}
We let $S = \sum_{k=1}^{q-1} \frac{ t^{k} \zeta^{k}}{k}$ and
regroup the terms, finding:
\begin{eqnarray*}
(1-t\overline \zeta^{q}) \cdot (\varphi(t) \cdot \zeta^{q} + S) & \equiv &
(1-t\zeta^{q}) \cdot (\varphi(t) \cdot \zeta^{q} + \overline   S) \mod q, \quad
\hbox{so} \\
-(\zeta^{q}-\overline \zeta^{q})\varphi(t) & \equiv
& (1-t\overline \zeta^{q}) S  - (1-t \zeta^{q}) \overline S  \mod q,
\end{eqnarray*}
and 
\[
-(\zeta^{q}-\overline \zeta^{q}) \varphi(t) \equiv \sum
\frac{t^{k}}{k}(\zeta^{k}-\overline \zeta^{k}) - \sum
\frac{t^{k+1}}{k }(\zeta^{k-q}-\overline \zeta^{k-q}) \mod q . \]
We regroup the powers of $\zeta$ using $q-k \equiv -k \mod q$, thus
$\zeta^{k-q}/k \equiv -\overline \zeta^{q-k}/(q-k)$, which can be
applied in the above for $k = 1, 2, \ldots, q-1$:
\begin{eqnarray*}
-(\zeta^{q}-\overline \zeta^{q}) \varphi(t)  \equiv 
 \sum_{k=1}^{q-1}  \ \frac{t^{k}-t^{2-k}}{k} \cdot 
(\zeta^{k}-\overline \zeta^{k}) \mod q,
\end{eqnarray*}
the statement of \rf{ponder}.
\end{proof}

Lemma \ref{tdev} yields essentially a system of equations modulo $q$
in the unknown $t$. It turns out that under some additional conditions
on $p$ and $q$, there are only three possible values for $t$ (one of
which is $t = 0$). The \textit{light } version of this condition was
presented in \cite{Mi1}; it reflects the main ideas which will
subsequently lead, by a more in depth study of the system \rf{ponder},
to a sharper inequality between $p$ and $q$. The light result is the
following:
\begin{proposition}
\label{qdivs}
Assume that $p > q$ are odd primes and there is a $\beta \in
\Q(\zeta)$ such that \rf{premq} holds.  Then
\begin{eqnarray}
\label{qdiv}
x+f \cdot y \equiv 0 \mod q^{2} \quad \hbox{ for some } \quad f \in
\{ -1, 0, 1 \} .
\end{eqnarray}
\end{proposition}
\begin{proof}
  Assume first that $x \equiv 0 \mod q$ and $x = q u$ with $(u, q)=1$.
  Since $(x,y) = 1$ and $p \neq$, it follows that $(x+\zeta^{a} y, q)
  = 1$, so the right hand side of \rf{premq} is a $q$ - adic integer.
  The equation is Galois - invariant, so we can replace $\zeta$ by
  $\zeta^{q}$.  Thus \rf{premq} becomes
\[ \frac{y+ q \overline \zeta^{q}  u}{y+q \zeta^{q} u} = \gamma^{q}, \]
with $\gamma = \pm \zeta^{2} \cdot \beta/\overline \beta$. Obviously
the above implies $\gamma \equiv 1 \mod q$, so $\gamma^{q} \equiv 1
\mod q^{2}$ and $y+q u \zeta^{2} \equiv y+ q u \overline \zeta^{2} \mod
q^{2}$, so $u \cdot (\zeta^{2}-\overline \zeta^{2}) \equiv 0 \mod q$. This
is only possible if $u \equiv 0 \mod q$ and thus $x \equiv 0 \mod
q^{2}$. Since we can interchange $x$ and $y$, this proves that if $x$
or $y$ is divisible by $q$, then it is divisible by $q^{2}$, which
takes care of $f = 0$ in this case. 

We may now assume that $x \cdot y \not \equiv 0 \mod q$ and use the
previous lemma, which implies that \rf{ponder} holds under the given
premises.  Since the set $\{\zeta, \zeta^{2}, \ldots, \zeta^{p-1}\}$
builds a base of the algebra $\Z[\zeta]/(q \cdot \Z[\zeta])$, the
coefficients of the single powers in the above identity must all
vanish and $p > q+1$ implies that the coefficient of $\zeta$ is $a_{1}
= t(1-t^{-4})$ and thus
\[ t^{4} \equiv 1 \mod q \] 
must hold. Furthermore, if $q+2 < p$, then the coefficient of
$\zeta^{2}$ is
\[ 2 \cdot a_{2} = (t^{2}-t^{-4}) \equiv 0 \quad \hbox{hence} \quad
t^{6} - 1 \equiv 0 \mod q. \] The last two congruences in $t$ have the
only common solution $t^{2} = 1 \mod q$. One easily verifies that if
this holds, then the right hand side in \rf{ponder} vanishes and thus
$\varphi(t) \equiv 0 \mod q$. This leads to the possible solution $x
\pm y \equiv 0 \mod q^{2}$; inserting the value back shows that this
is indeed a solution of \rf{premq}. If $p = q+2$, then we still have
$a_{1} = t^{-3}(t^{4}-1)$ so $t^{4} \equiv 1 \mod q$. If $t^{2} - 1
\equiv 0 \mod q$, we find the previous solution. So let us assume that
$t^{2} \equiv -1 \mod q$ and consider the second coefficient: but
$\varphi(t) \overline \zeta^{q} = \varphi(t) \zeta^{2}$ has in this
case a contribution to $a_{2}$. We estimate this coefficient by using
$t^{2} \equiv -1 \mod q$:
\begin{eqnarray*}
 2 \cdot a_{2} & \equiv & t^{2}-t^{-4} - 2 \varphi(t) \equiv  
-t^{-4}\left(t^{6} - t^{2}+t^{2}-1+2 t^{4} \varphi(t)\right) \\ 
& \equiv & 
t^{2} - 1 + 2\varphi(t) \equiv 2(\varphi(t) -1)\mod q , 
\end{eqnarray*}
a congruence which is satisfied by $\varphi(t) \equiv 1 \mod q$. We
have to consider also 
\begin{eqnarray*}
3 \cdot a_{3} & = & (t^{3} - t^{-5}) - (t^{q-1} - t^{-q-1}) \equiv 0
\mod q \quad \Leftrightarrow \\
0 & \equiv & t^{-5}(t^{8}-1) - (1 - t^{-2}) \mod q.
\end{eqnarray*}
But if $t^{2} \equiv -1 \mod q$, then the first term vanishes while the
second is $-2 \not \equiv 0 \mod q$, so $t^{2} \equiv -1 \mod q$ is
not possible. This takes care also of the case $p = q+2$, thus
completing the proof of the proposition.
\end{proof}

It follows from Corollary \ref{crhodef} that
\begin{corollary}
\label{c0}
  If $p > q > 3$ are odd primes for which \rf{FC} has non trivial
  solutions and such that $q \nmid h_{p}^{-}$, then \rf{qdiv} holds.
\end{corollary}
\begin{proof}
The premises of Corollary \ref{crhodef} are given and thus \rf{rhod}
holds. By setting $\beta = \rho_{1}$ in this equation, we find that
the hypotheses of Proposition \ref{qdivs} also hold, and by its proof it
follows that \rf{qdiv} must be true.
\end{proof}
\subsection{Sharpening}
Let $\rg{k}$ be a field and $\id{T}$ be the space of sequences on
$\rg{k}(t)$. We define the following operators on $\id{T}$:
\begin{eqnarray}
\label{thetadef}
b_{n} & = & \theta_{+}(a_{n}) = a_{n} - t \cdot a_{n-1} \nonumber \\
c_{n} & = & \theta_{-}(a_{n}) = t \cdot a_{n} - a_{n-1} \nonumber \\
d_{n} & = & \Theta(a_{n}) = \theta_{+} (\theta_{-}(a_{n})). \nonumber
\end{eqnarray}
Furthermore we let $\Delta$ be, classically, the forward difference
operator $\Delta \  a_{n} \ = \ a_{n} - a_{n-1}$ and $n^{\underline{k}} = n
\cdot (n-1) \ldots (n-k+1)$ be the \nth{k} falling power of $n$, so
$\Delta n^{\underline{k}} = k \cdot (n-1)^{\underline{k-1}}$. With,
the main properties of the operators in \rf{thetadef} are given by
\begin{lemma}
\label{thetaprop}
The operators $\theta_{+}, \theta_{-}$ are linear and they commute,
thus $\Theta = \theta_{+} \circ \theta_{-} = \theta_{-} \circ
\theta_{+}$. Furthermore, 
\begin{eqnarray}
\label{op1}
\begin{array}{c c c c c c c}
\theta_{+}(t^{n}) & = & 0 & \hbox{ and } & \theta_{+}(t^{-n}) & = &
(1-t^{2}) t^{-n}, \\ 
\theta_{-}(t^{-n}) & = & 0 & \hbox{ and } & \theta_{-}(t^{n}) & = &
-(1-t^{2}) t^{n-1},
\end{array}
\\
\label{op2}
\begin{array}{c c l}
\theta_{+}^{l}(n^{\underline{k}} \cdot t^{n}) & = & \frac{k!}{l!}
\cdot (n-l)^{\underline{k-l}} \cdot t^{n}, \\ 
\theta_{-}^{l}(n^{\underline{k}} \cdot t^{-n}) & = & \frac{k!}{l!} \cdot
(n-l)^{\underline{k-l}} \cdot t^{-(n-l)},
\end{array}
\end{eqnarray} 
where we set $a^{\underline{k-l}} = 0$ if $k < l$. In particular, we
have:
\begin{eqnarray}
\label{op3}
\begin{array}{c c l}
\theta_{+}^{k}(n^{\underline{k}} \cdot t^{n}) & = & k! \cdot t^{n}, \\ 
\theta_{-}^{k}(n^{\underline{k}} \cdot t^{-n}) & = & k! \cdot
t^{-(n-k)}, \\
\Theta^{k}(n^{\underline{k}} \cdot t^{n}) & = & k! \cdot (t^{2}-1)^{k}
\cdot t^{n-k}, \\
\Theta^{k}(n^{\underline{k}} \cdot t^{-n}) & = & k! \cdot (-1)^{k}
\cdot (t^{2}-1)^{k} \cdot t^{-(n-k)}.
\end{array}
\end{eqnarray} 
\end{lemma}
\begin{proof}
Commutativity follows by a straight forward computation from 
\[ \theta_{+} \circ \theta_{-} (a_{n} ) = \theta_{-} \circ \theta_{+}
(a_{n} ) = t \cdot (a_{n}+a_{n-2}) - (t^{2}+1) a_{n-1} .\] The rules
\rf{op1} are also easily verified and they yield \rf{op2} by induction
on $k$. Finally, the first two actions in \rf{op3} are obtained by
setting $l = k$ in \rf{op2}, while the action of $\Theta$ is obtained
due to commutativity, by setting $\Theta^{k} = \theta_{-}^{k} \circ
\theta^{k}_{+}$ or $\Theta^{k} = \theta_{+}^{k} \circ
\theta^{k}_{-}$, depending whether the operand is $t^{n}$ or
$t^{-n}$. Note that $k+1$ consecutive values of $a_{n}$ are necessary
for applying $\theta_{\pm}^{k}$, while $\Theta^{k}$ requires $2k+1$
consecutive values.
\end{proof}
The task we pursue is to improve our estimates on pairs $p, q$ for
which the system \rf{ponder} has no other solutions except \rf{qdiv};
in particular, we are concerned with $p < q$ - since Proposition
\ref{qdivs} deals already with $p > q$. We shall use the fact on which
the proof of Proposition \ref{qdivs} relays: $(\zeta^{k})_{k=1}^{p-1}$
form a base of the algebra $\Z[\zeta]/(q \Z[\zeta])$ and this allows
one consider \rf{ponder} as a linear system modulo $q$. Concretely,
the coefficients of $\zeta^{k} - \overline \zeta^{k}$ in that equation
must vanish, for $k = 1, 2, \ldots, \frac{p-1}{2}$. Let $0 < \nu <
\frac{p-1}{2}$ be the value for which $\nu \equiv q \mod p$ or $\nu
\equiv -q \mod p$; then, with $\delta_{ij}$ the Kronecker $\delta$,
the above remark yields the equations:
\begin{eqnarray}
\label{lsys}
-\delta_{\nu,k} \cdot \varphi(t) & \equiv &
\sum_{j \geq 0; jp + k < q} \frac{t^{k+pj} - t^{2-(k+pj)}}{pj+k} \\ &
- & \sum_{j \geq 0; jp + (p-k) < q} \frac{t^{p-k+pj} - 
t^{2-(p-k+pj)}}{p-k+jp} \mod q. \nonumber
\end{eqnarray}
The index value $\nu$ is singular for the equations above; first, it
is the only index for which the equations are not homogeneous. Second
the number of terms in the sums of the right hand side changes between
$0 <k < \nu$ and $p/2 > k > \nu$. In these two intervals \rf{lsys}
yields homogeneous equations which manifest in the vanishing of
polynomials of fixed degree in $k$. This suggests the use of the
difference operators defined above. Let $5 \leq p < q$ be primes. We
shall take the approach of choosing the one of the intervals $0 < k <
\nu$ or $\nu < k < p/2$, which has more elements: in these intervals
\rf{lsys} translates into polynomial equations of the type $f_{q}(k;
t) = 0$. Having a contiguous interval on which this equation holds,
one can use the iteration of $\Theta$ in order to reduce the degree in
$k$ of the polynomial $f_{q}$. We have thus to distinguish the cases
$\nu < p/4$ and $\nu > p/4$ \footnote{\ One may also take the approach
  of considering the whole interval $0 < k < p/2$; in this case the
  polynomials $f_{q}(k; t)$ change the degree and shape when $k$
  passes the "singular" value $k = \nu$. The computations become more
  intricate, for a gain of a factor at most $2$. We choose to analyze
  here the simpler approach.}.
\begin{proposition}
\label{p+}
Let $5 \leq p < q$ be primes such that \rf{ponder} holds and $\nu$ be
defined above. Suppose that $\nu > p /4$; if additionally, $q <
\frac{p^{2}}{16}$, then \rf{qdiv} holds. 
\end{proposition}
\begin{proof}
Let $n = \left[q/p\right]$. The equation \rf{lsys} yields on the
interval $0 < k < \nu$: 
\begin{eqnarray*}
\sum_{0 \leq j \leq n} \frac{t^{k+pj} - t^{2-(k+pj)}}{pj+k}
\equiv \sum_{0 \leq j < n} \frac{t^{p-k+pj} - 
t^{2-(p-k+pj)}}{p-k+jp} \mod q. 
\end{eqnarray*}
After eliminating denominators, this yields a polynomial equation:
\begin{eqnarray*}
(-1)^{n} k^{2n} & \cdot & \sum_{0 \leq j \leq n} \left(t^{k+pj} -
  t^{2-(k+pj)} 
\right) + O(k^{2n-1})  \equiv \\ (-1)^{n-1} k^{2n} & \cdot & \sum_{0
  \leq j < n} \left(  t^{p-k+pj} - t^{2-(p-k+pj)} \right) + O(k^{2n-1}).
\end{eqnarray*}
In order to eliminate the lower order terms in $k$, we may take
$\Theta^{2n}$ on both sides of the congruence. This requires at least
$2(2n)+1$ contiguous points, so $1 \leq k -2n < k+2n < p/4$, which
means $2(2n)+1 < p/4$. If this is provided, the equation reduces,
after simplifying by $(-1)^{n} \cdot (2n)! \cdot (1-t^{2})^{2n}$, to :
\begin{eqnarray}
\label{redgen1}
& & \sum_{0 \leq j \leq n} \left(t^{k+pj-2n} - t^{2-(k+pj-2n)}\right)
\\ & + & 
\sum_{0 \leq j < n} \left(  t^{p-k+2n+pj} - t^{2-(p+2n-k+pj)} \right)
  \equiv 0 \mod q. \nonumber
 \end{eqnarray}
If $t \not \in \{-1, 0, 1\}$ then we can apply $\theta_{+}$ and
$\theta_{-}$ independently to the above congruence. This yields:
\begin{eqnarray*}
0 & \equiv & \sum_{0 \leq j \leq n} t^{2-(k+pj-2n)} - \sum_{0 \leq j < n}
t^{p-k+2n+pj} \quad \hbox{and } \\
0 & \equiv & \sum_{0 \leq j \leq n} t^{k+pj-2n}-  
\sum_{0 \leq j < n}t^{2-(p+2n-k+pj)},
\end{eqnarray*} 
and, upon multiplication by the lowest power of $t$, 
\begin{eqnarray}
\label{toadd}
0 & \equiv & \sum_{0 \leq j \leq n} t^{pj} - \sum_{0 \leq j < n}
t^{p(n+1)-2+pj} \mod q \quad \hbox{and } \\
0 & \equiv & \sum_{0 \leq j \leq n} t^{pn+pj-2}- \sum_{0 \leq j < n}  
t^{pj} \mod q .
\nonumber 
\end{eqnarray}
Adding up the two congruences, we obtain $t^{pn} \equiv -t^{pn-2}
\mod q$ with the solutions $t \equiv 0$ and $t^{2} \equiv -1 \mod q$.
We show that the latter solution is impossible by reinserting it in
\rf{redgen1}; this yields, after simple computations, $t^{k}+t^{-k}
\equiv 0 \mod q$. Since we assumed $t \not \equiv 0$, it follows that
$(-1)^{k}+1 \equiv 0 \mod q$. It suffices to take $k$ even in order to
reach a contradiction. Let us finally examine all the conditions on
$\nu$ (and thus on $p$ and $q$), which allowed us to reach this
contradiction. Adding the points necessary for the final application
of $\theta_{\pm}$ together with the condition that $k$ be even, we
find:
\[ 2n+1 \leq k \leq p/4-(2n+1) ,\]
condition which is satisfied by the even value $k = 2(n+1)$, provided
that $4n+3 < p/4$. On the other hand, we find from the definition of
$\nu$ and the fact that $\nu > p/4$, that $p(4n+3) > 4q$, and thus
\[ p^{2}/4 > p(4n+3) > 4q, \] as claimed.
\end{proof}
\begin{proposition}
\label{p-}
Let $5 \leq p < q$ be primes such that \rf{ponder} holds and $\nu$ be
defined above. Suppose that $\nu < p /4$; if additionally, $q <
\frac{p(p-20)}{16}$, then \rf{qdiv} holds.
\end{proposition}
\begin{proof}
The proof of this proposition follows the same line as the previous,
but raises few particular obstructions. We shall let 
\[ n = \begin{cases}
\lfloor q/p \rfloor & \hbox{if} \quad (q \mod p) < p/4, \\
\lfloor q/p \rfloor + 1 & \hbox{if} \quad (q \mod p) > 3 p/4. 
\end{cases}
\]
 The equation \rf{lsys} yields now on the
interval $ \nu < k < p/4$: 
\begin{eqnarray*}
\sum_{0 \leq j \leq n} \frac{t^{k+pj} - t^{2-(k+pj)}}{pj+k}
\equiv \sum_{0 \leq j \leq n} \frac{t^{p-k+pj} - 
t^{2-(p-k+pj)}}{p-k+jp} \mod q. 
\end{eqnarray*}
Note that there are equally many terms in the sums of both sides of
the above congruences, unlike the case of the previous proposition.
This perpetuates down to the analog of \rf{toadd}, in which the two
congruences become identical; they both yield the condition
\begin{eqnarray}
\label{condcum}
 t^{p(n+1)} \equiv 1 \mod q \quad \hbox{or} \quad t^{p(n+1)} \equiv
 t^{2} \mod q,
\end{eqnarray}
whose deduction is left to the reader. Note that this condition is
equivalent to applying any of $\theta_{+} \Theta^{2n+1}$ or
$\theta_{-} \Theta^{2n+1}$ to the original system \rf{lsys}.

In order to draw a contradiction we shall  have to consider lower
order terms in $k$. Let
\[ \sigma_{j} = t^{k+pj} - t^{2-(k+pj)} \quad \hbox{ and } \quad \tau_{j} = 
t^{p-k+pj} -t^{2-(p-k+pj)}; \]
with some additional work, the first congruence
yields, after elimination of denominators:
\begin{eqnarray*}
\sum_{0 \leq j \leq n} \sigma_{j} \cdot
  \left(k^{\underline{2n+1}}-\left[(n+j+1)p -(2n+1)n\right] \cdot
  k^{\underline{2n}}\right)  & + & \\ 
\sum_{0 \leq j \leq n} \tau_{j} \cdot
  \left(k^{\underline{2n+1}}-\left[(n-j)p -(2n+1)n\right] \cdot
  k^{\underline{2n}}\right) & + & O(k^{2n-1})  \equiv  0 \mod q.
\end{eqnarray*}
We apply $\Theta^{2n}$ to the above and let 
\[ \sigma'_j = t^{k-2n+pj} - t^{2-(k-2n+pj)} \quad \hbox{ and } \quad
\tau'_{j} = t^{p-k+2n+pj} -t^{2-2n-(p-k+pj)}. \]
With this we obtain
\begin{eqnarray}
\label{theta2n}
\sum_{j=0}^{n} & & \sigma'_{j} \left((2n+1) k - \left[(n+j+1)p
    -(2n+1)n\right] \right) + \\  \sum_{j=0}^{n}  & & \tau'_{j}
\left((2n+1) k - \left[(n-j)p  -(2n+1)n\right] \right) \equiv 0 \mod
q. \nonumber
\end{eqnarray}
We now apply $\theta_{+}^{2}$ to the above relation; this cancels the
terms in $t^{k}$ and modifies the terms in $t^{-k}$. Note that by
commutativity, $\theta^{-} \theta_{+}^{2} \Theta^{2n} = \theta_{+}
\Theta^{2n+1}$, which yields \rf{condcum}. But applying $\theta^{-}$
after $\theta_{+}^{2}$ to \rf{theta2n} yields to a cancellation of all
but the terms in $\theta_{+}^{2} (k t^{-k})$; conversely, it is
precisely these terms which are canceled if the condition
\rf{condcum} holds.  Since $\theta_{+}^{2} t^{-k} \not \equiv 0 \mod
q$ if $t \mod q \not \in \{-1,0,1\}$, this eventually leads to the
congruence:
\begin{eqnarray*}
\begin{array}{l c c c c}
p \cdot t^{2n} \cdot \left(\sum_{j=0}^{n} \left((j+1)t^{2-pj} + j
  t^{p+pj}\right) 
  \theta_{+}^{2}(t^{-k}) \right) & \equiv & 0 & \mod q, & \hbox{so} \\
\sum_{j=0}^{n} \left((j+1)t^{2-pj} + j  t^{p+pj}\right)  & \equiv & 0 &
  \mod q, & \hbox{and} \\
t^{2-pn} (n+1) \cdot \sum_{j=0}^{n} t^{pj} + t^{p} \cdot
  (1-t^{2-p(n+1)}) \cdot \sum_{j=0}^{n} j t^{pj} & \equiv & 0 & \mod
  q. &  
\end{array}
\end{eqnarray*}
We can now reintroduce the alternative \rf{condcum} in the last
congruence above. If $t^{p(n+1)} \equiv 1 \mod q$, then the first sum
$\sum_{j=0}^{n} t^{pj}$ vanishes; furthermore, since $(n,q) = 1$, one
easily verifies that $\sum_{j=0}^{n} t^{pj} \equiv \sum_{j=0}^{n} j
t^{pj} \equiv 0 \mod q$ cannot simultaneously hold, and thus it
follows that $t^{2-p(n+1)} \equiv 1$. Since we also assumed
$t^{(n+1)p} \equiv 1 \mod q$, it follows that $t^{2} \equiv 1 \mod q$,
as required. Suppose now that in \rf{condcum} it is the condition
$t^{p(n+1)} \equiv t^{2} \mod q$ which holds; by inserting this in the
last congruence above, we find (since $t(n+1) \not \equiv 0 \mod q$)
that $\sum_{j=0}^{n} t^{pj} \equiv 0 \mod q$ and we are in the
previous case. Both ways, it follows that $t \mod q \in \{-1, 0, 1\}$.

We finally have to derive the inequality between $p$ and $q$, for
which the proof above holds. The condition is that the interval $(p/4,
p/2)$ contains sufficient contiguous points for applying both
$\theta_{\pm} \Theta^{2n+1}$ and $\theta_{+}^{2} \Theta^{2n}$; i.e.
$4n+5 < p/4$. Note that by definition of $n$, we always have $np > q$
and thus the previous inequality amounts to $4q < 4np < p(p/4-5)$ and
thus
\[ \frac{p(p-20)}{16} > q. \]
This completes the proof of the proposition.  
\end{proof}
\subsection{Proof of Theorem \ref{main}}
The statement of Theorem \ref{main} follows directly from Corollary
\ref{c0} together with the sharpening Propositions \ref{p+} and
\ref{p-}.

\subsection{The Resulting Case Analysis}
We suppose that the Fermat - Catalan equation \rf{FC} has a solution
for odd primes $p, q$ with $p \not \equiv 1 \mod q$, $q \nmid
h_{p}^{-}$ and $\max \{p, \frac{p(p-20)}{16}\} > q$. Then Theorem
\ref{main} holds and we are reduced to investigate the case $e = 0$
(Case I) or $e = 1$ (Case II) each with three subcases: $f = -1$ (case
a), $f = 0$ (case b) and $f = 1$ (case c); together, this yields the
Table 1, with six cases.
 \begin{table}
\caption{Cases of Fermat - Catalan}
\vspace{2pt}
\begin{tabular} {|c|c|c|c|}
\hline
& {\bf a ( f = -1)}&
{\bf b (f = 0)} &
{\bf c (f = 1)}
 \\[2pt]
\hline
{\bf I (e = 0)}   &  I a     &  I b &     I c \\
{\bf II (e = 1)}   &  II a     &  II b &     II c \\[2pt]
\hline
\end{tabular}
\end{table}
Furthermore, if either $q \nmid h(p, q)$ or $q \nmid h_{pq}^-$, Corollary
\ref{crhodef} holds and in particular the identity \rf{rhod}. We aim
next to eliminate the unit $\varepsilon$ in this identity, using the
fact that, by Proposition \ref{lemmaA}, the $q$ - primary units of
$\Q(\zeta)$ are global \nth{q} powers. This shall be done by a case by
case study. In view of Lemma \ref{lemmaC}, we let $0 < m < q$ with
$m(p-1) \equiv 1 \mod q$ and $\gamma = \left(
  \frac{(1-\zeta)^{p-1}}{p}\right)^{m}$ be the unit in \rf{uc}.

Suppose first that $e = 0$ and thus $\alpha = x+y \zeta$. In case a,
$x \equiv y \mod q^{2}$ and
\[ \varepsilon \cdot \rho^{q} = \alpha \equiv x (1+\zeta) \mod q^{2}
.\] But then $\delta = \frac{\varepsilon}{1+\zeta} =_{q} x$ and by
Lemma \ref{unitp} it follows that $\delta \in \Z[\zeta]^{q}$;
consequently $\alpha = (1+\zeta) \rho^{q}$ in this case. If $f = 0$
(case b)), then $x \equiv 0 \mod q^{2}$ and $\varepsilon =_{q} y$ and
Lemma \ref{unitp} shows that $\varepsilon$ is a \nth{q} power, so
$\alpha = \rho^{q}$ in this case. Finally, if $f = 1$, then
$\varepsilon =_{q} -y(1-\zeta)$. Since $\delta = \varepsilon : \gamma
=_{q} y$, the Lemma \ref{lemmaC} implies that $\delta$ is a \nth{q}
power. It follows that $\alpha = \gamma \cdot \rho^{q}$ in this case.

We assume next that $e = 1$ and consider the three possible values of
$f$. In case a, $\alpha = \frac{x+\zeta y}{1-\zeta} \equiv
\frac{x(1+\zeta)}{1-\zeta}$ and combining the Lemmata \ref{unitp} and
\ref{lemmaC} we find that $\alpha = (1+\zeta) / \gamma
\rho^{q}$. Likewise, $\alpha = \rho^{q}/\gamma$ in case b) and $\alpha
= \rho^{q}$ in case c. 
 \begin{table}
\caption{Values of the unit $\varepsilon$ in the six Cases}
\vspace{2pt}
\begin{tabular} {|c|c|c|c|}
\hline
& {\bf a }&
{\bf b } &
{\bf c }
 \\[2pt]
\hline
{\bf I }   &  $(1+\zeta)$      &  $1$ &     $\gamma$ \\
{\bf II }   &  $\frac{1+\zeta}{\gamma} $     &  $1/\gamma$ &      $1$ \\[2pt]
\hline
\end{tabular}
\end{table}
We combine all these results in Table 2 and the following
\begin{proposition}
\label{propqpow}
Let $p, q$ be odd primes with $p \not \equiv 1 \mod q$, $q \nmid h(p,
q)$ and suppose that $\max \{p, \frac{p(p-20)}{16}\} > q$ and the
Fermat - Catalan equation \rf{FC} has a non trivial solution.
Furthermore, let $0 < m < q$ be an integer with $m(p-1) \equiv 1 \mod
p$ and $\gamma = \left( \frac{(1-\zeta)^{p-1}}{p}\right)^{m} \in
\Z[\zeta]^{\times}$. Then for $e \in \{0, 1\}$ and $f \in \{-1, 0,
1\}$ like in Table 1, the following identity holds (with
$\delta_{a,b}$ being the Kronecker $\delta$ symbol):
\begin{eqnarray}
\label{sixqpow}
\alpha = \frac{x+\zeta y}{(1-\zeta)^{e}} =
(1+\zeta)^{\delta_{f,-1}} \cdot \gamma^{\delta_{f,1}-e} \cdot
\rho^{q},  \quad \hbox{with} \quad \rho \in \Q(\zeta). 
\label{qpowcomb}
\end{eqnarray}
\end{proposition} 

We proceed with a case by case analysis of possible solutions in the
above six cases. The results come in different levels of complexity
and require different class number conditions - essentially the two
possibilities $q \nmid h(p,q)$ or $q | h_{pq}^-$, mentioned above. 

The simplest fact is that in three out of six cases, one deduces some
Wieferich - type local conditions, involving only the exponents $p$
and $q$. This is the topic of the next section. In the following
section, keeping the same class number condition, we show that one can
give lower bounds on $\max\{|x|, |y|\}$: this is Theorem \ref{main1}.
We sharpen subsequently the class number condition to $q \nmid
h_{pq}^{-}$ and prove, by a generalization of Kummer descent - as used
by Kummer in his Theorem on the Second Case of Fermat's Last Theorem,
\cite{Wa} - that two additional cases ($f = -1$) are impossible. This
leaves on last case - which we called the \textit{Ast\'erisque} - Case
- untreated by conditions involving only the exponents $p, q$. By
using the sharper class number condition, we are able to improve the
lower bound in this case to one on the \textit{minimum} $\min\{|x|,
|y|\}$; this is Theorem \ref{main2}, which is a first generalization
of Catalan's conjecture. Finally, by applying this Theorem together
with an additional consequence of the Kummer descent, we prove the
Theorem \ref{trc} on the rational case of Catalan's equation, which is
the most exhaustive result of this paper, since it shows the lack of
solutions of Catalan's equation in the rationals, provided some
conditions hold, which are related only to the exponents.

\section{The Wieferich Cases}
We assume in this section that \rf{FC} has non-trivial solutions for
odd prime exponents $p, q$ for which the premises of Theorem
\ref{main} hold. Based on this theorem, we can thus assume that the
solutions are in one of the cases given in the above tables. The three
simplest cases lead to some Wieferich - type (see \cite{Ri})
condition.
\begin{proposition}
\label{p1}
Notations being as above, if $e = 0$ and $f = -1$, then 
\[ 2^{q-1} \equiv 1 \mod q^2. \]
Furthermore, $X =_{q} Y =_{q} 1$.
\end{proposition}
\begin{proof}
  Since $e = 0$, we are in Case I and $X+Y = A^q$; also, $f = -1$
  means $X-Y \equiv 0 \mod q^2$, so $X =_q Y$ and $X+Y =_q 2X =_q A^q
  =_q 1$. But from \rf{FC}, $X^p + Y^p = Z^q =_q 2X^p =_q 1$. Dividing
  the last two relations, we find $X^{p-1} =_q 1$ and since $q \nmid
  p-1$ by hypothesis, it follows also that $X =_q 1$. Combined with
  $2X =_q = 1$ this yields the statement of the proposition. Since $X
  =_{q} Y$ by definition of this case and $2X =_{q} 2 =_{q} 1$, the
  second statement follows too.
\end{proof}
\begin{proposition}
\label{p2}
Notations being as above, if $e = 1$ and $f = 0$, then 
\[ p^{q-1} \equiv 1 \mod q^2 \]
and $Y =_{q} 1$.
\end{proposition}
\begin{proof}
Since $e = 1$, we are in Case II and $\frac{X^p+Y^p}{X+Y} = p \cdot
B^q =_q p$; also, $f = 0$ means $X \equiv 0 \mod q^2$, so 
\[ \frac{X^p+Y^p}{X+Y} =_q Y^{p-1} =_q p. \]
But $X^p +Y^p = Z^q =_q Y^p =_q 1$ and since $(p,q) = 1$, we must have
$Y =_q 1$, which is the second statement of the Proposition. Combined
with the previous equivalence, this yields $Y =_q p =_q 1$, which
leads to the first claim.
\end{proof}
The third Wieferich case has a more complex statement. This is:
\begin{proposition}
\label{p3}
Notations being as above, if $e = 1$ and $f = -1$, then 
\[ \left(2^{p-1} \cdot p^{p}\right)^{q-1} \equiv 1 \mod q^2 \]
and $X =_{q} Y =_{q} p^{m}$, where $m (p-1) \equiv 1 \mod q$.
\end{proposition}
\begin{proof}
  Since $e = 1$, we are still in Case II, so $\frac{X^p+Y^p}{X+Y} = p
  \cdot B^q =_q p$; also, $f = -1$ implies $X =_{q} Y$ and
  $\frac{X^p+Y^p}{X+Y}=_{q} X^{p-1} =_{q} p$, the second claim of the
  Proposition. Furthermore, $X^{p}+Y^{p} =_{q} 2X^{p} =_{q} z^{q}
  =_{q} 1$. Raising this to the power $p-1$ and the preceding
  equivalence to the power $p$, we find after division that:
\[ 2^{p-1} X^{p(p-1)} =_{q} 2^{p-1} p^{p} =_{q} 1. \]
This is the first statement of the Proposition and completes the
proof.
\end{proof}

\section{Lower Bounds and Proof of Theorem \ref{main1}}
We assume in this section that \rf{FC} has non trivial solutions for
odd primes $p, q$ with $p \not \equiv 1 \mod q$, $q <
\max\{\frac{p(p-20)}{16}\}$ and such that $q \nmid h(p, q)$. The
purpose of this section is to prove Theorem \ref{main1}.

The following $q$ - adic expansion will serve for gaining estimates in
all six cases under investigation.
\begin{lemma}
\label{vanish}
Let $\rho \in \id{O}(\KL)^{\times}$ be an algebraic integer with $q$ -
adic expansion
\[ \rho = a \cdot \sum_{m=0}^{\infty} \binom{1/q}{m} (\mu b)^{m},
\quad \hbox{with} \quad \mu \in \{\zeta, 1/(1\pm \zeta)\}, \quad a \in
\Z_{q}, \ b \in \Q^{\times}, \ v_{q}(b) = k \geq 2. \] Furthermore,
suppose there is a real number such that $|\sigma(\rho)| \geq M > 0$
for all $\sigma \in G_{p}$. Then
\begin{eqnarray}
\label{lbgen}
M \geq \frac{1}{p-1} \cdot
\frac{q^{(p-2)\left(k-\frac{q}{q-1}\right)}}{2^{p-2} } .
\end{eqnarray}
\end{lemma}
\begin{proof}
Let 
\[ \nu = (\zeta^{2} - \zeta) \cdot \begin{cases}
1 & \hbox{if} \quad \mu = \zeta \\
\mu^{-(p-3)} & \hbox{otherwise}.
\end{cases}
\]
It is an easy verification, that $\Tr_{\K/\Q} \left(\nu \cdot
  \mu^{i}\right) = 0$ for $j = 0, 1, \ldots, p-3$. The case $\mu =
\zeta$ is trivial, since $\Tr\left(\zeta^{i+2}-\zeta^{i+1}\right) =
(-1)-(-1) = 0$ for $0 \leq i \leq p-2$. If $\mu = 1/(1 \pm \zeta)$,
then
\begin{eqnarray*}
 \Tr\left(\nu \cdot \mu^{j}\right) & = & \Tr\left((\zeta^{2}-\zeta) \cdot
  (1\pm\zeta)^{p-3-j}\right) \\ & = & \sum_{i=0}^{p-3-j} (\pm 1)^{i}
   \binom{p-3-j}{i} \cdot \Tr\left((\zeta^{2}-\zeta) \zeta^{i}\right)
   = 0, 
\end{eqnarray*}
as claimed. Let now $\delta = \nu \cdot \rho \in \id{O}(\K)^{\times}$.
The $q$ - adic expansion of $\rho$ together with the above remark on
the trace of $\nu \cdot \mu^{j}$ shows that the first $p-2$ terms in
the $q$ - adic expansion of $\Delta = \Tr_{\K/\Q}(\delta) \in \Z$
vanish. Thus 
\begin{eqnarray}
\label{Delta}
\Delta = \binom{1/q}{p-2} b^{p-2} \cdot \left(\Tr
  \left(\mu^{p-2}\right)+O(q)\right).
\end{eqnarray} 
Note that 
\begin{eqnarray}
\label{binomvq}
v_{q}\left(\binom{1/q}{n}\right) = v_{q}\left(\frac{1}{q^{n}} \cdot
  \frac{(1-q) \cdot \ldots \cdot (1-(n-1)q)}{n!}\right) = -n-v_{q}(n!),
\end{eqnarray}
and since $v_{q}(b) \geq k$, it follows that 
\begin{eqnarray*}
 v_{q}\left(\binom{1/q}{p-2} b^{p-2} \right) & \geq & k(p-2) - (p-2) -
v_{q} \left((p-2)!\right) \\ & > & (p-2)(k-1-1/(q-1)) =
(p-2)\left(k-\frac{q}{q-1}\right) .
\end{eqnarray*}
We now show that $\Delta \neq 0$. Assume first that $\mu =
\zeta$. Then from \rf{Delta} we find that $\Delta = \binom{1/q}{p-2}
b^{p-2} \cdot (p + o(q))$. Since $b \neq 0$ the first $q$ - adic term
is non vanishing, and so $\Delta \neq 0$. The proof is similar for
$\mu = 1 \pm \zeta$. Finally, $\Delta$ is a rational integer and by
assembling all the information, we find:
\[ q^{(p-2)\left(k-\frac{q}{q-1}\right)} \leq |\Delta| \leq
\sum_{\sigma \in G_{p}} | \sigma(\nu \cdot \rho) | \leq 2^{p-2} \cdot
(p-1) \cdot M. \]
The claim follows from these inequalities.
\end{proof}

The proof of Theorem \ref{main1} follows now from the Lemma and the
fact that $x+f y \equiv 0 \mod q^{2}$.
\begin{proof}
We shall show in a case by case analysis, that a $q$ - adic expansion
such as required by Lemma \ref{vanish} exists. If $f = 1$, then
\rf{sixqpow} and \rf{uc} yield: 
\[ \frac{x+\zeta y}{(1-\zeta)^{e}} =_{q} \frac{1+\zeta}{\gamma^{e} } =
p^{(m-q)e} \frac{1+\zeta}{(1-\zeta)^{e} } \cdot \rho^{q} , \]
and thus $\rho \in \Z[\zeta]$ with 
\[ \rho^{q} = p^{(q-m)e} \cdot y \left(1 + \frac{x-y}{y(1+\zeta)}\right)
. \] If $e = 0$, the expansion follows by Proposition \ref{p1}, since
the leading term is $y =_{q} 1$; if $e = 1$ the leading term is
$p^{(q-m)} y =_{q} 1$, by Proposition \ref{p3}. We still have to
deduce the bound $M$ from the expression of $\rho$. But
\begin{eqnarray*}
 \left| \rho \right| & = & \left|  p^{(q-m)e} \cdot \frac{x+\zeta
    y}{1+\zeta} \right|^{1/q} \leq p^{(1-m/q) e} \cdot
\left(\frac{|x|+|y|}{1/p}\right)^{1/q} \\ 
& \leq & \left(2p^{1+(q-m) e} \cdot \max\{|x|,|y|\}\right)^{1/q} .
\end{eqnarray*} 
the last estimate is obviously Galois invariant, so we can replace $M$
in the Lemma \ref{vanish} by this value. It follows that
\[ M = \left(2p^{1+(q-m) e}
  \cdot \max\{|x|,|y|\}\right)^{1/q} \geq \frac{1}{p-1} \cdot
\frac{q^{(p-2)\left(\frac{q-2}{q-1}\right)}}{2^{p-2} } , \] and
\begin{eqnarray*}
%\label{comest}
  \max\{|x|,|y|\} \geq \frac{1}{2} \cdot \left(\frac{1}{p(p-1)} \cdot
\left(\frac{q^{\frac{q-2}{q-1}}}{2 }\right)^{p-2}\right)^{q}
\end{eqnarray*}
for both cases, as claimed in the first inequality of \rf{lowb}. 

Let now $f = 0$ and $y \equiv 0 \mod q$, to fix the ideas. Then
$\frac{x+\zeta y}{(1-\zeta)^{e}} = \gamma^{-e} \rho$ and $\rho^{q} =
p^{(q-m)e} (x+\zeta y)$. If $e = 1$, the leading term is $p^{(q-m)e} x
=_{q} 1$ by Proposition \ref{p2}, otherwise, $A^{q} = x+y =_{q} x
=_{q} 1$, so the expansion of $\rho$ follows in both cases.
Furthermore,
\[ \left| \rho \right| = \left|  p^{(q-m)e} \cdot (x+\zeta
  y) \right|^{1/q} \leq \left(p^{(q-m) e} \cdot 2 \cdot \max\{|x|,
  |y|\}\right)^{1/q} = M .\] Like before, by applying the Lemma
\ref{vanish}, we find that \rf{lowb} holds.

Finally, if $f = 1$, some usual computations yield
\[ \rho^{q} = \frac{x+\zeta  y}{1-\zeta} \cdot p^{me} = -y \cdot
p^{me} \cdot \left(1 - \frac{x+y}{1-\zeta}\right). \] It follows
immediately from the fact that the parenthesis on the right hand side
of the last identity is a $q$ - adic \nth{q} power (since $x+y \equiv
0 \mod q^{2}$) that so must then be the cofactor $yp^{me}$; this shows
the existence of the $q$ - adic expansion of $\rho$ required by Lemma
\ref{vanish}. The details for the estimation of $M$ are analogous to
the case $f = -1$ and are left to the reader.
\end{proof}

\section{Kummer Descent}
We shall prove in this section the following main Theorem, which
generalizes Kummer's descent method to the present context.
\begin{theorem}
\label{tdesc}
Let $p, q > 3$ be primes such that $-1 \in <p \mod q>$, $\zeta,
\xi \in \C$ be respectively \nth{p} and \nth{q} primitive roots of
unity, $\KL = \Q(\zeta, \xi)$ and $\KL^{++}$  the fixed field of the
partial complex conjugations $\jmath_p, \jmath_q$. Suppose that the
equation
\begin{equation}
\label{desc}
  X^{q} + Y^{q} = \varepsilon \cdot \lambda^{N} \cdot {\lambda'}^{M}  
\cdot Z^{q}
\end{equation}
admits solutions with $X, Y, Z \in \id{O}\left( \KL^{++}\right), X
\cdot Y \cdot Z \neq 0$ and $(X \cdot Y \cdot Z, p \cdot q) = (1)$ and
$X, Y$ are not units. Here $\lambda = (\xi-\overline \xi), \lambda' =
(\zeta - \overline \zeta)$ and $\varepsilon \in \id{O}\left(
\KL^{++}\right)^{\times}$; $M, N$ are integers with $N > 2q, N$ even
and $M = 0$ or $M \geq 2$. Then $Z$ is not a unit and $q | h_{pq}^-$.
\end{theorem}

The next Lemma will explain the condition $-1 \in <p \mod q>$:
\begin{lemma}
\label{realp}
Let $p, q$ be odd primes and $\K' = \Q(\xi)$ be the \nth{q} cyclotomic
extension. Then $p$ splits in $\K'$ in real prime ideals iff $-1 \in
<p \mod q>$.
\end{lemma}
\begin{proof}
  This is a direct consequence of Kummer's Theorem on the splitting of
  primes in extensions with a power base for the ring of algebraic
  integers \cite{La}. Let $\Phi_{q}(X) = \prod_{i=1}^{q-1}
  (X-\xi^{i})$ be the \nth{q} cyclotomic polynomial, $n = \ord_{q}(p)
  = |< p \mod q>| $ and $F(X) = \prod_{j=0}^{n-1} (X-\xi^{p^j}) \in
  \Z[\xi][X]$. If $\rg{k} = \F_{p^{n}}$ is the finite field with
  $p^{n}$ elements, then $\rg{k}$ is the smallest field of
  characteristic $p$ which contains a non trivial \nth{q} root of
  unity. Let $\rho \in \rg{k}$ be such a root of unity. Then there is
  a natural map $\iota : \id{O}(\K') \rightarrow \rg{k}$ given by $\xi
  \mapsto \rho$. Let then $\tilde{f}(X) = \iota(F(X)) \in \F_{p}(X)$
  and $f \in \Z[X]$ be some polynomial with $\tilde{f} = f \mod p$.
  Then $\tilde{f} \in \F_{p}[X]$ is an irreducible factor of $\Phi(X)
  \mod p $; if $\eu{p} = (f(\xi), p)$, then $\eu{p}$ is a prime above
  $(p)$ and each prime above $(p)$ arises in this way, by a choice of
  $\rho \in \rg{k}$. In particular, $\xi \mod \eu{p} = \rho$ and
  $\iota$ is in fact the reduction $\mod \eu{p}$ map. Furthermore,
  $\tilde{f} = \iota\left(F(X)\right) = f(X) \mod p$, where in general
  only the second polynomial has rational integer coefficients.
  
  After this exposition of Kummer's Theorem, we can proceed with the
  proof of our Lemma. First note that since $\K'/\Q$ is a CM Galois
  extension, all the primes above $(p)$ are simultaneously real or not
  real. Let us first suppose that $\eu{p} = (f(\xi), p)$ is a real
  ideal. Since $\tilde{f}(X) = F(X) \mod \eu{p}$ and $\overline
  {\eu{p}} = \eu{p}$, it follows that $F(X) = \overline {F(X)}$ and
  under the action of $\iota$,
\[ \tilde{f}(X) = \prod_{i=0}^{n-1} \left(X-\rho^{p^{i}}\right) =   
\prod_{i=0}^{n-1} \left(X-\rho^{-p^{i}}\right)\] But $\rho \in
\rg{k}$, which is a field in which $\tilde{f}(X)$ has unique
decomposition. Thus $\rho^{-1} \in \{ \rho^{p^{i}} : i = 0, 1, 2
\ldots, n-1\} $ and by the definition of $n$ it follows that $-1 \in <
p\mod q>$, as claimed. Conversely, if $-1 \in <p \mod q>$, then $F(X)
= \overline {F(X)}$ and it follows that $\eu{p}$ is invariant under
complex conjugation.
\end{proof}
We proceed with the proof of the Theorem, assume that $q \nmid
h_{pq}^{-}$ holds under the given hypotheses and will derive a
contradiction.  The two statements of the theorem are apparently
contradictory: if we show that $q \nmid h_{pq}^-$ is impossible, then
it is irrelevant whether $Z$ is a unit or not.  For technical reasons,
however, it will be useful to show that under the given premises, if
$X, Y$ are not units, then neither is $Z$. Note that by Proposition
\ref{lemmaA}, it follows that $E_{q} = E^{q}$ and $\id{C} =
\id{C}^{q}$, with $\id{C}$ the ideal class group of $\KL$. The quite
lengthy proof is a straightforward adaption of the descent method used
by Kummer in the proof of his fundamental Theorem on the Second Case
of Fermat's Last Theorem (see \cite{Wa}, Chapter 9). The additional
problems are linked to the fact that we work in a larger field.
  
We start with a simple fact:
\begin{lemma}
\label{snorm}
Let $X, Y \in \KL^{++}$ verify the premises of the theorem; in
particular, suppose that $X^q+Y^q \equiv 0 \mod \lambda^N \cdot
{\lambda'}^M$ and $v_{\lambda}(X^{q}+Y^{q}) = N,
v_{\lambda'}(X^{q}+Y^{q}) = M$. Then $v_{\lambda}(X+Y) = N-(q-1)$ and
$v_{\lambda'}(X+Y) = M$.
\end{lemma}  
\begin{proof}
  Any integer $\gamma \in \id{O}(\KL)$ has the $\lambda$ - development
\[ \gamma = \sum_{i=0}^G g_i \cdot \lambda^i, \quad \hbox{for some } 
\quad G \in \N, \] and 
\[ g_i = \sum_{j=1}^{p-1} g_{i,j}
\zeta^j\in \Z[\zeta], \quad 0 \leq g_{i,j} < q .\] Let $X = x_0 + x_1
\cdot \lambda + O(\lambda^2)$. Since $\jmath_q(\lambda) = -\lambda$
and $X \in \KL^{++}$, so $\jmath_q(X) = X$, we must have $x_1 = 0$, so
$X = x_0 + O(\lambda^2)$. Likewise, $Y = y_0 + O(\lambda^2)$. Thus
\rf{desc} implies that $x_0^q+y_0^q \equiv 0 \mod \lambda^2$ and since
$(X \cdot Y, q) = 1$, we have $(x_0/y_0)^q \equiv -1 \mod
\lambda^2$. Since $\Z[\zeta]/(q \Z[\zeta])$ contains no \nth{q} roots
of unity except $1$, it follows that $x_0/y_0 \equiv -1 \mod
\lambda^2$ and $x_0 + y_0 \equiv X+Y \equiv 0 \mod \lambda^2$. The
algebraic integers 
\begin{eqnarray}
\label{phis}
\phi'_{i} & = & \xi^{i} X + \overline \xi^{i}
  Y \in \KL, \quad \hbox{for} \quad i=1, 2, \ldots, q-1
\end{eqnarray}
have the common divisor 
\[ (\phi'_{i}, \phi'_{j}) = \left((\xi^{i}-\xi^{j}) \cdot Y,
  (\overline \xi^{i}-\overline \xi^{j}) \cdot X\right) = (\lambda) .
\] But $(X+Y) \cdot \prod_{i=1}^{q-1} \phi'_{i} = X^{q}+Y^{q} \equiv 0
\mod \lambda^{N}$. Thus $v_{\lambda}(X+Y) = N-(q-1)$, as claimed.

Due to $(\phi'_{i}, \phi'_{j}) = (\lambda)$, it follows also that if
$\eu{P} | (\lambda')$ is a prime ideal of $\KL^{++}$ with $\eu{P} |
(X+Y)$, then $\eu{P}^{M} | (X+Y)$. The primes above $(p)$ in
$\eu{\KL}$ are $(\eu{p}, (1-\zeta))$ for some prime $\eu{p} \in \K'$
and the hypothesis together with Lemma \ref{realp} imply that they are
real primes. Suppose that $M > 0$ (there is nothing to prove for $M =
0$!); then for $\eu{P} | (p)$ we have $X^{q}+Y^{q} \equiv 0 \mod
\eu{P}$ and $(-X/Y)^{q} \equiv 1 \mod \eu{P}$. Since $\eu{P}$ is a
prime ideal, $\id{O}(\KL^{++})/\eu{P}$ is a field and there is an
integer $0 \leq a < q$ such that $-X/Y \equiv \xi^{a} \mod \eu{P}$.
Taking complex conjugates - under consideration of the fact that $X/Y$
is invariant under conjugation - we also have $-X/Y \equiv \xi^{-a}
\mod \overline {\eu{P}}$. But since $\eu{P} = \overline {\eu{P}}$, it
follows that $\xi^{2a} \equiv 1 \mod \eu{P}$ and $a = 0$. This holds
for all primes above $(p)$ and together with the previous remark
implies that $v_{\lambda'}(X+Y) = M$, as claimed.
\end{proof}

We wish to normalize the algebraic integers defined in \rf{phis},
eliminating all primes above $q$ and $p$. Using the result of the
Lemma \ref{snorm}, this can be done as follows:
\begin{eqnarray}
\label{phis1}
\phi_{i} & = & \frac{\xi^{i} X + \overline \xi^{i}
  Y}{\xi^{i}-\overline \xi^{i}} , \quad \hbox{for} \quad i=1, 2,
  \ldots, q-1, \\ \phi_{0} & = & \frac{q(X+Y)}{\varepsilon \cdot
  \lambda^{N} \cdot {\lambda'}^M}.
\end{eqnarray}
It follows from the Lemma \ref{snorm} that $(\phi_{i}, p \cdot q) =
(1)$ for $i=0, 1, \ldots, q-1$ and
\begin{eqnarray*}
\prod_{i=0}^{q-1} \phi_{i} & = & \frac{1}{q} \cdot \phi_{0} \cdot
\prod_{i=1}^{q-1} \phi'_{i} = \frac{q}{q \cdot \varepsilon \cdot
\lambda^{N}\cdot {\lambda'}^M} \cdot (X+Y) \cdot \prod_{i=1}^{q-1}
\phi'_{i} = \frac{X^{q}+Y^{q}}{\varepsilon \cdot \lambda^{N}\cdot
{\lambda'}^M}.
\end{eqnarray*}
Finally, this yields
\begin{eqnarray}
\label{prod2}
\prod_{i=0}^{q-1} \ \phi_{i} = Z^{q}.
\end{eqnarray}
If $\KL_p \subset \KL$ is the subfield fixed by $\jmath_p$, the
definition of $\phi_i$ implies $\phi_i \in \KL_p$ for $i > 0$ and
\rf{prod2} shows that this holds also for $\phi_0$:
\begin{eqnarray}
\label{inkp}
\phi_i \in \KL_p \quad i = 0, 1, \ldots, q-1.
\end{eqnarray} 
From $(\phi'_{i}, \phi'_{j}) = (\lambda)$ and the definition of
$\phi_{i}$ we deduce that
\begin{eqnarray}
\label{coprim}
(\phi_{i}, \phi_{j}) = 1 \quad \hbox{for} \quad i, j \in \{0, 1, \ldots,
q-1\} , \quad \hbox{and} \quad i \neq j.
\end{eqnarray}

We want to show that $\phi_i$ are not units. This implies that $Z$ is
not a unit, as a consequence of \rf{prod2}. For $\phi_0$, this fact
will follow indirectly, with more work. We prove it first only for $i
> 0$ and investigate the $q$ - expansion of $\phi_i$:
\begin{eqnarray}
\label{modq}
 \phi_{i} & = & \frac{\xi^{i} X + \xi^{-i} Y}{\xi^{i}-\xi^{-i}} =
\frac{\xi^{i}(X+Y)-(\xi^{i}-\xi^{-i})\cdot Y}{\xi^{i}-\xi^{-i}} \\ & =
& -Y +  \frac{\xi^{i}(X+Y)}{\xi^{i}-\xi^{-i}}   =_{q}  -Y \nonumber.
\end{eqnarray}
Note that $\phi_{q-i} = -\overline \phi_i$. If $\phi_i$ is a unit, then
$\delta = \phi_i/\phi_{q-i} = \phi_i/\overline \phi_i$ is a root of
unity and by \rf{modq}, since $Y \in \R$, it follows that $\delta =_q
1$. By Lemma \ref{lru} it follows that $\delta = \pm \zeta^a$ for some
$a \in \ZM{p}$. Then $\zeta^{-a/2} \phi_i = \pm \zeta^{a/2}
\phi_{q-i}$ and a short computation shows that
\[ X = - Y \cdot \frac{\overline \zeta^{a/2} \cdot \overline \xi \mp  
\zeta^{a/2} \cdot 
\xi}{\overline \zeta^{a/2} \cdot \xi \mp \zeta^{a/2} \cdot \overline
\xi} = \gamma \cdot Y. \] But $\gamma$ is a unit and $X = \gamma Y$ is
a contradiction to $(X, Y) = (1)$, since $X$ and $Y$ are not
units. The contradiction confirms our claim that $\phi_i$ are not
units, for $i > 0$; thus $Z$ is not a unit.

We can now apply Lemma \ref{idealq} with $n = q, C = Z, \KL' = \KL_p$,
thus obtaining:
\begin{lemma}
\label{l2}
Let the premises of Theorem \ref{tdesc} hold, the normalized elements
$\phi_{i}, i > 0$ be defined by \rf{phis1} and the ideals $\eu{A}_{i}
= (\phi_{i}, Z)$. If $\ \KL_p \subset \Q(\zeta, \xi)$ is the subfield
fixed by $\jmath_p$, then $\eu{A}_{i}$ are principal and there are
$\mu'_{i} \in \id{O}(\KL_p)$ and $\eta_{i} \in
\left(\id{O}(\KL_p)\right)^{\times}$, such that \rf{alg} holds.
\end{lemma}
Note that if $\phi_0$ is not a unit, the result of Lemma \ref{l2} also
holds for $\phi_0$. We proceed our proof, allowing for both
possibilities. It will turn out that the same computations which allow
descent also imply in the long run that $\phi_0$ is not a unit.

From \rf{prod2}, $\phi_{0} = \frac{Z^{q}}{\prod_{i=1}^{q-1} \phi_{i}}
=_{q}(-Y)^{1-q} =_{q} -Y$. If $\phi_0$ is not a unit, we saw that we
can write $\phi_0 = \eta_0 \cdot \mu_0^q$, with $\eta_i, \mu_i$ like
in \rf{alg}; otherwise we may set $\phi_0 = \eta_0$. In both cases,
the unit $\eta_0$ is defined and $\eta_0 =_q \phi_0 =_q -Y$.

Thus, for all $i=0, 1, \ldots, q-1$, we have $\eta_i =_q -Y$; the
Lemma \ref{unitp} implies that $\eta_{i}$ must be \nth{q} powers, so
\begin{eqnarray}
\label{alg1}
\phi_{i} & = & \mu_{i}^{q}, \quad \hbox{for} \quad i = 0,
\ldots, q-1.
\end{eqnarray}
If $\phi_0$ is a unit, then the previous remarks imply that $\phi_{0}
= \mu_{0}^{q}$, with $\mu_{0}$ a unit of the same field. Otherwise, by
the same reasoning as in the case $i > 0$, $\mu_0 \in
\id{O}\left(\KL_{p}\right)$.

We are prepared for the main computations which will allow to perform
the descent. We evaluate $\phi_{i} \times \phi_{-i}$ for $i > 0$,
using the identity in \rf{modq}:
\begin{eqnarray*}
\psi_{i} = \phi_{i} \times \phi_{-i} & = & \left(-Y +
  \frac{\xi^{i}(X+Y)}{\xi^{i}-\xi^{-i}} \right) \cdot  \left(-Y +
  \frac{\xi^{-i}(X+Y)}{\xi^{-i}-\xi^{i}} \right) \\
& = &\left(\frac{X+Y}{1-\overline{\xi}^{2i}}-Y\right) \cdot  
\left(\frac{X+Y}{1-\xi^{2i}}-Y\right) \\
& = & Y^{2} + \left(\frac{X+Y}{|1-\xi^{2i}|}\right)^{2} -Y\cdot(X+Y)
  \cdot
  \left(\frac{1}{1-\overline{\xi}^{2i}}+\frac{1}{1-\xi^{2i}}\right) \\
& = &  \left(\frac{X+Y}{|1-\xi^{2i}|}\right)^{2} - X\cdot Y.
\end{eqnarray*}
The last equation above shows that $\psi_i \in \id{O}(\KL^{++})$. We
let $\psi_0 = \phi_0^2$, so we also have $\psi_0 \in \id{O}(\KL^{++})$:
\begin{eqnarray}
\label{real}
\psi_i \in \id{O}(\KL^{++}) \quad \hbox{for} \quad i = 0, 1, \ldots, (q-1)/2.
\end{eqnarray} 
By subtracting the values of $\psi$ for two indices $i \not \equiv \pm
j \mod q$ we find $\psi_{i}-\psi_{j} = \delta_{i,j} \cdot (X+Y)^{2}$,
with $\delta_{i, j} = 1/|(1-\xi^{2i}| - 1/|(1-\xi^{2j}|$. For the
choice of such indices we need here that $q \geq 5$. We claim that
$\lambda^{2} \cdot \delta_{i,j} = \eta_{i,j} \in
\id{O}\left(\KL^{++}\right)^{\times}$. Indeed,
\begin{eqnarray}
\label{deltaij}
\lambda^{2} \cdot \delta_{i,j} & = &
  \frac{\lambda^{2}}{|(1-\xi^{2i})(1-\xi^{2j})|^{2}} \cdot
  \left(|1-\xi^{2j}|^{2} - |1-\xi^{2i}|^{2}\right) \nonumber \\ 
& = & \frac{\lambda^{2}}{|(1-\xi^{2i})(1-\xi^{2j})|^{2}} \cdot
  \left((2-\xi^{2i}-\xi^{-2i})-(2-\xi^{2j}-\xi^{-2j})\right) \\
& = & \frac{\lambda^{2} \cdot (\xi^{2j}-\xi^{2i}) \cdot (1-\overline
  \xi^{2(i+j)})}{|(1-\xi^{2i})(1-\xi^{2j})|^{2}} \nonumber
\end{eqnarray}
In our definition $\lambda = \xi - \overline \xi$ is an imaginary
number, so $\lambda^{2}$ is real and so is $\lambda^{2} \cdot
\delta_{i,j}$. The last equality above shows that
$v_{\eu{q}}\left(\lambda^{2} \cdot \delta_{i,j}\right) = 0$ and since
it is real and invariant under $\jmath_{q}$ it follows that
$\eta_{i,j} \in \id{O}(\KL^{++})^{\times}$, as claimed.

We now substitute the definition \rf{phis} of $\phi_{i}$ and \rf{alg1}
in the recent results, finding:
\begin{eqnarray*}
\psi_{i} - \psi_{j} & = & \eta^{2} \cdot \left((\mu_{i} \cdot
  \mu_{q-i})^{q} - (\mu_{j} \cdot \mu_{q-j})^{q}\right) \\ & = &
  \eta_{i,j} \cdot \lambda^{-2} \cdot (X+Y)^{2} = \eta_{i,j} \cdot
  \lambda^{-2} \cdot \left(\frac{\varepsilon \cdot \lambda^{N} \cdot
  {\lambda'}^M \cdot \phi_{0}}{q}\right)^{2} \\ & = & \eta_{i,j}
  \cdot \lambda^{-2} \cdot \left(\frac{\varepsilon \cdot \lambda^{N}
  \cdot {\lambda'}^M \cdot \eta \cdot \mu_{0}^{q}}{q}\right)^{2} \\ & = &
  \left(\eta_{i,j} \cdot \eta^{2} \cdot
  \left(\frac{ \lambda^{q-1}}{q}\right)^{2}\right) \times
  {\lambda'}^{2M} \cdot\lambda^{2(N-q)} \cdot \mu_{0}^{2q}.
\end{eqnarray*}
After division by $\eta^{2}$ this yields:
\begin{eqnarray}
\label{down}
(\mu_{i} \cdot \mu_{q-i})^{q} - (\mu_{j} \cdot  \mu_{q-j})^{q}  = 
\eta' \cdot {\lambda'}^{M'} \cdot \lambda^{N'} \cdot \mu_{0}^{2q},
\end{eqnarray}
where $\eta' = \left(\delta_{i,j} \cdot
  \left(\frac{\lambda^{q-1}}{q}\right)^{2}\right) \in
\id{O}\left(\KL^{++}\right)^{\times}$ and $N' = 2(N-q) = N + (N-2q) >
N$ is even, $M' = 2M$. Also, by \rf{real}, the numbers occurring at
the \nth{q} power in \rf{down} are elements of $\KL^{++}$.
  
We have shown that for $i > 0$, $\phi_i$ are not units, and thus the
\nth{q} powers on the left hand side of \rf{down} are neither units.
We still have to show that $\mu_0$ is not a unit. For this we write
$X' = \mu_{i} \cdot \mu_{q-i}, Y' = -\mu_{j} \cdot \mu_{q-j} \not \in
\Z[\zeta, \xi]^{\times}$ and $Z' = \mu_0^2$ and anticipate the next
descent step. We start from \rf{down}, which can be rephrased to
\begin{eqnarray*}
{X'}^{q} + {Y'}^{q}  = \eta' \cdot {\lambda'}^{M'} \cdot \lambda^{N'}
\cdot {Z'}^{q}.
\end{eqnarray*}
Since $X', Y'$ are not units, we have proved above that it follows
that $Z'$ is not a unit either, so $\mu_0, \phi_0$ are indeed not
units, as claimed.

In order to control the descent, we let $I$ be the group of integer
ideals of $\Z[\zeta,\xi]$ and consider $\omega : I \rightarrow \N$,
the function counting the number of distinct prime ideals which divide
an ideal $\eu{A} \in I$. Thus $\omega$ generalizes the analogous
function defined on the integers. We show that $0 <
\omega(\mu_{0}^{2}) = \omega(\mu_{0}) < \omega(Z)$. Indeed, since the
$\phi_{i}$ are coprime and not units, relation \rf{prod2} together
with \rf{alg1} imply that
\[ \omega(Z) = \sum_{i=0}^{q-1} \omega(\phi_{i}) = \sum_{i=0}^{q-1}
\omega(\mu_{i}) ,\] thus $\omega(\mu_{0}) < \omega(Z)$. The inequality
$\omega(\mu_{0}) > 0$ rephrases the fact the $\mu_0$ is not a unit,
which we have proved. We thus have the main argument of descent:
\begin{proposition}
\label{pdesc}
Let $p, q; X, Y, Z; \varepsilon, \lambda, \lambda', N, M$ be like in
the statement of Theorem \ref{tdesc} and suppose that $q \nmid
h_{pq}^-$. Then there are $X^{(1)}, Y^{(1)}, Z^{(1)} \in
\id{O}(\KL^{++})$, a unit $\varepsilon^{(1)} \in
\id{O}\left(\KL^{++}\right)^{\times}$, an even integer $N^{(1)} >
N$ and $M^{(1)} \geq M$, such that
\begin{eqnarray}
\label{desca}
  \left(X^{(1)}\right)^{q} + \left(Y^{(1)}\right)^{q} =
  \varepsilon^{(1)} \cdot \lambda^{N^{(1)}} \cdot {\lambda'}^{M^{(1)}}
  \cdot \left(Z^{(1)}\right)^{q}
\end{eqnarray}
and $Z^{(1)} \mid Z, X^{(1)}, Y^{(1)}, Z^{(1)} \not \in \Z[\zeta,
\xi]^{\times}$. Finally, $\omega\left(Z^{(1)}\right) < \omega(Z)$,
where $\omega$ is the distinct prime factor counting function.
\end{proposition}
\begin{proof}
  With the notations above, we let $X^{(1)} = \mu_i \cdot \mu_{q-i}$,
  $Y^{(1)} = -\mu_j \cdot \mu_{q-j}$ and $Z^{(1)} = \mu_0^2$; also
  $\varepsilon^{(1)} = \eta'$ and $N^{(1)} = N'$. We have proved that
  $X^{(1)}, Y^{(1)}, Z^{(1)} \in \id{O}(\KL^{++})$, they are coprime
  an non vanishing and $(X^{(1)} \cdot Y^{(1)} \cdot Z^{(1)}, p \cdot
  q) = (\phi_0, \cdot q) = (1)$. Also, $\eta' \in
  \id{O}\left(\KL^{++}\right)^{\times}$ and $N' = 2(N - q) > N > 2q$
  is even; $M' = 2M \geq M$, trivially. It was shown that $X^{(1)},
  Y^{(1)}, Z^{(1)} \not \in \Z[\zeta, \xi]^{\times}$. Thus all the
  conditions of Theorem \ref{tdesc} are verified. The equation
  \rf{desca} is then a reformulation of \rf{down}.
\end{proof}
The proof of Theorem \ref{tdesc} follows now easily:
\begin{proof}
The Proposition \ref{pdesc} can be applied recursively to \rf{desca},
thus generating an infinite sequence 
\[ (Z) = (Z^{(0)}) \subset (Z^{(1)}) \subset (Z^{(2)}) \subset \ldots \subset 
(Z^{(k)}) \subset \ldots , \] such that $\omega(Z^{(k)}) >
\omega(Z^{(k+1)})$ for all $k \geq 0$. But $Z=Z^{(0)}$ has only a
finite number of prime factors and the function $\omega$ is positive
integer valued, so it cannot decrease indefinitely. This is a
contradiction which shows that the hypothesis $q \nmid h_{pq}^-$ of
Proposition \ref{pdesc} is untenable, thus proving Theorem
\ref{tdesc}.
\end{proof}

\section{Case Analysis}
We consider two primes $p, q > 3$ such that $q \nmid h_{pq}^-$ and $-1
\in <p \mod q>$ and suppose that \rf{FC} holds for these values of $p,
q$. The Barlow - Abel relations imply then that
\[  \frac{X^p+Y^p}{X+Y} = p^e \cdot A^q, \]
for some $e \in \{0,1\}$ and $A \in \Z$. By theorem \ref{main}, 
\begin{eqnarray*}
x+f \cdot y & \equiv & 0 \mod q^{2} \quad \hbox{ with } \quad f \in
\{-1,0,1\}. 
\end{eqnarray*}
Together this yields six cases, three of which have been dealt with
above, by means of Wieferich relations. We shall investigate below the
remaining cases.

\subsection{The Descent Cases}
\begin{theorem}
\label{desc1}
Notations being as above and assuming the premises of Theorem
\ref{main2}, the equation \rf{FC} has no solution with $e=1, f = -1$.
\end{theorem}
\begin{proof}
  Assume that \rf{FC} has a solution with $e=1, f=-1$. Then $(X+Y)/p$
  is a \nth{q} power, divisible by $p \cdot q$. Let $v_q(X+Y) = nq$
  and $v_p(X+Y) = mq-1$, so
\[ X+Y = p^{mq-1} \cdot q^{nq} \cdot C^q \quad C \in \Z, \ (C, p q) = 1. \]
By \rf{sixqpow} we have in the present case:
\begin{eqnarray}
\label{rholaq}
 \frac{x+\zeta y}{1-\zeta} = -y + \frac{x+y}{1-\zeta} = \rho^{q} .
\end{eqnarray}
Note that 
\begin{eqnarray*}
 \alpha - \overline \alpha &= & \frac{x+y}{1-\zeta} -
\frac{x+y}{1-\overline \zeta} = \frac{(x+y)(1+\zeta)}{1-\zeta} = -
\frac{(x+y)(\zeta^a+\overline \zeta^a)}{\zeta^a-\overline \zeta^a}
\quad \hbox{so} \\ \rho^{q} - \overline \rho^{q} & = &
\prod_{j=0}^{q-1} (\xi^{j} \cdot \rho -\xi^{-j} \cdot \overline \rho)
= - C^{q} \cdot (p^m \cdot q^n)^q \cdot \frac{\zeta^a+\overline
\zeta^a}{p(\zeta^a-\overline \zeta^a)},
\end{eqnarray*}
with $a =(p+1)/2$.

We define the following system of normed divisors of $C^q$:
\begin{eqnarray}
\label{phisA}
\phi_i & = & \frac{\xi^i \rho - \overline \xi^i \overline \rho}{\xi^i
- \overline \xi^i} \\ 
\label{phi0A}
\phi_0 & = & -\frac{\rho-\overline \rho}{p^{mq-1}
\cdot q^{nq-1}} \cdot \frac{\zeta^a-\overline
\zeta^a}{\zeta^a+\overline \zeta^a}.
\end{eqnarray}
Let 
\begin{eqnarray*}
\varepsilon_1 & = & \prod_{c=1}^{p-1} \frac{\zeta^c - \overline
\zeta^c}{\lambda'} \in \Z[\zeta+\overline \zeta]^{\times} \subset
\id{O}\left(\KL^{++}\right)^{\times} , \\ \varepsilon_2 & = &
\prod_{c=1}^{q-1} \frac{\xi^c - \overline \xi^c}{\lambda} \in
\Z[\xi+\overline \xi]^{\times} \subset
\id{O}\left(\KL^{++}\right)^{\times}.
\end{eqnarray*}
Then $p = \varepsilon_1 \cdot {\lambda'}^{p-1}$ and $q = \varepsilon_2
\cdot {\lambda}^{q-1}$ and \rf{phi0A} can be rewritten as
\[\rho-\overline \rho = \varepsilon \cdot \lambda^{(q-1)(nq-1)} \cdot 
{\lambda'}^{(p-1)(mq-1)-1} \cdot \phi_0, \]
with 
\[\varepsilon = \varepsilon_1^{mq-1} \cdot \varepsilon_2^{nq-1} \cdot 
\frac{(\zeta^a+\overline \zeta^a) \lambda'}{\zeta^a-\overline \zeta^a}
\in \id{O}\left(\KL^{++}\right)^{\times} .\]
Finally, with $N = (q-1)(nq-1)$ and $M = (p-1)(mq-1)-1$, we have 
\begin{eqnarray}
\label{phi0B}
\rho-\overline \rho = \varepsilon \cdot \lambda^N \cdot {\lambda'}^M
\cdot \phi_0 \quad \hbox{with} \quad \varepsilon \in
\id{O}\left(\KL^{++}\right)^{\times}.
\end{eqnarray}
The rest of the proof goes through a series of steps which were proved
in detail in the previous section, so we list the arguments, leaving
it to the reader to check the details.

We have by construction $\prod_{i=0}^{q-1} \phi_i = C^q$ and since
$(C, pq) = 1$, a fortiori $(\phi_i, pq) = (1)$. Since $(\rho,
\overline \rho) = (1)$, one verifies that $(\phi_i, \phi_j) = (1)$ for
$0 \leq i \neq j < q$. We can apply Lemma \ref{idealq} to $\phi_i$,
with $m = q, \KL' = \KL_p$, and find
\[  \phi_i = \eta_i \cdot \mu_i^q  \quad \eta_i \in \id{O}(\KL_p)^{\times},
 \ \mu_i \in \id{O}(\KL_p). \]
We have from \rf{rholaq} that $\rho^q =_q -y$; then there is an
integer $t$ with $t^q \equiv -y \mod q^{nq}$ and the $q$-adic
development based on \rf{rholaq} yields $\rho \equiv \overline \rho
\equiv t \mod q^{nq-1}$. But for $i > 0$ we have $\eta_i =_q \phi_i
=_q t$ and since $\phi_0 = C^q/\prod_{i > 0} \phi_i$, we also have
$\eta_0 =_q \phi_0 =_q t$. Consequently, we may assume that
\[ \phi_i = \eta_0 \mu_i^q, \]
and $\eta_0 \in \Z[\zeta+\overline \zeta]^{\times} \subset
\id{O}(\KL^{++})^{\times}$.

Finally we define $\psi_i = \phi_i \cdot \phi_{q-i} \in
\id{O}(\KL^{++})$ and verify that for $i \not \equiv \pm j \mod q$ we have
\begin{eqnarray*}
 \psi_i + \psi_j & = & \eta_0^2 \cdot
 \left(\left(\mu_i\cdot\mu_{q-i}\right)^q + \left(\mu_j \cdot
 \mu_{q-j}\right)^q\right) = \eta_{i,j} \cdot \lambda^{-2} \cdot (\rho
 - \overline \rho)^2 \\ & = & \eta_{i,j} \cdot \left(
 \frac{\varepsilon \cdot \lambda^N \cdot {\lambda'}^M}{ \lambda }
 \right)^2 \times \left(\eta_0^2 \cdot \mu_0^{2q} \right),
\end{eqnarray*}
where $\eta_{i,j} = \lambda^2 \cdot \delta_{i,j}$ is the unit in
\rf{deltaij}. After dividing by $\eta_0^2$, we set $X = \mu_i \cdot
\mu_{q-i}, Y = \mu_j \cdot \mu_{q-j}, Z = \mu_0^2, N' = 2(N-1), M' = 2M,
\varepsilon' = \varepsilon^2 \cdot \eta_{i,j}$ and find:
\[X^q + Y^q = \varepsilon' \cdot \lambda^{N'} \cdot {\lambda'}^{M'} 
\cdot Z^q. \] The hypotheses $(X,Y,Z) = (XYZ,pq) = (1)$, $X, Y, Z,
\varepsilon' \in \KL^{++}$, $N > 2 q$ is even and $M \geq 0$ being all
fulfilled, as has been showed above, we can apply the Kummer descent
Theorem \ref{tdesc}. This raises a contradiction to $q \nmid
h_{pq}^-$, which proves the statement of this Proposition.
\end{proof}

Next we treat the case $p \nmid z, q | x+y$:
\begin{theorem}
\label{desc2}
Notations being as above and assuming the premises of Theorem
\ref{main2}, the equation \rf{FC} has no solution with $e=0, f = 1$.
\end{theorem}
\begin{proof}
This is a case with $(X+Y,p) = (1)$ and $M = 0$ in the descent
theorem. By Corollary \ref{crhodef}, we have $\alpha = \varepsilon
\cdot \rho^q$. The $q$ - adic development of $\rho$ is more delicate
in this case and we shall work it out in detail - the rest of the
proof being exempt of surprises. We are in the First Case and 
\begin{eqnarray}
\label{xpy}
x+y \equiv 0\mod q^q,
\end{eqnarray} so
\[\frac{x^p+y^p}{x+y} =_q p y^{p-1} = B^q =_q 1. \]
If $m(p-1) = 1 + nq$, \rf{sixqpow} yields in this case, for a $\rho$
twisted by a root of unity:
\begin{eqnarray*}
 \rho^q  = - \zeta^{-m/2} \cdot (1-\zeta)^{-nq} \cdot y p^m \cdot \left(1 -
\frac{x+y}{(1-\zeta) y p^m}\right).
\end{eqnarray*}
Since the cofactor of $y p^{m}$ is a $q$ - adic \nth{q} power, it
follows from the above equation that $y p^m =_q 1$, so there is a $t
\in \Z$ with $t^q \equiv -y p^m \mod q^{2q}$. The ring $\ZM{q}[\zeta]$
contains no non trivial \nth{q} roots of unity (since $q$ in not
ramified in $\Z_{q}[\zeta]$), so the resulting $q$ - adic extension of
$\rho$ starts as follows:
\[ \rho = \frac{t}{\zeta^{m/2q} (1-\zeta)^m} \cdot \left(1 - 
\frac{x+y}{(1-\zeta) q y p^m} + O(q^{2(q-1)}) \right). \] From the
definition it follows that $n$ is odd and one verifies that $\rho$
verifies the necessary condition $(\rho/\overline\rho)^q \equiv -\zeta
\mod q^{q}$.

We now investigate an adequate factoring of $x+y = B^q$. We have
\[ (\zeta^{-1/2} + \zeta^{1/2}) (x+y) = \zeta^{-1/2} \alpha + \zeta^{1/2} 
\overline \alpha = \varepsilon \cdot \left(\zeta^{-1/2} \rho^q +
\zeta^{1/2} \overline \rho^q\right). \] Defining $\rho_1 = \overline
\zeta^{1/2q} \cdot \rho$, we have $\rho_1/\overline \rho_1 \equiv -1
\mod q^{q-1}$ and
\[   \varepsilon \cdot (\rho_1^q + \overline \rho_1^q) = (\zeta^{1/2}+
\overline \zeta^{1/2}) \cdot B^q: \]
this looks like a good starting point. Let $q^{nq} || (x+y)$, $C =
B/q^n$ with $(C, pq) = 1$ and define
\begin{eqnarray*}
\phi_i & = & \frac{\xi^i \rho_1 + \overline \xi^i \overline \rho_1
}{\xi^i-\overline \xi^i} \quad \hbox{ for } \quad i = 1, 2, \ldots,
q-1 \quad \hbox{ and } \\ \phi_0 & = & \frac{\varepsilon (\rho_1 + \overline
\rho_1)}{q^{nq-1} \cdot (\zeta^{1/2}+\overline \zeta^{1/2})}.
\end{eqnarray*}
Note that $\phi_0$ is an algebraic integer, since $\rho_1/\overline
\rho_1 \equiv -1 \mod q^{nq-1}$. Then
\begin{eqnarray*} 
\prod_{i=0}^{q-1} \phi_i  =  \left(\frac{\rho_1^q + \overline
\rho_1^q}{q}\right) \times \left(\frac{\varepsilon}{q^{nq-1} \cdot
(\zeta^{1/2}+\overline \zeta^{1/2})}\right) = (B/q^n)^q = C^q
\end{eqnarray*}
According to the usual frame, one verifies that $(\phi_i, \phi_j) =
(1)$ for $0 \leq i \neq j < q$ and by Lemma \ref{idealq} it follows that 
\[ \phi_i = \eta_i \cdot \mu_i^q, \quad \hbox{ with } \quad \eta_i \in 
\id{O}(\KL_p)^{\times}, \ \mu_i \in \id{O}(\KL_p).  \] Since $\phi_i
\equiv -\overline \rho_1 \mod \left(\frac{\rho_1 + \overline
\rho_1}{\lambda}\right)$ and $\phi_0 = C^q/\prod_{i > 0} \phi_i$, it
follows that $\phi_i =_q \eta_i =_q -\rho_1$ and, using Lemma
\ref{lru}, one deduces, after eventual modification of $\mu_i$, that
\[ \phi_i = \eta_0 \cdot \mu_i^q \quad i = 0, 1, \ldots, q-1. \]
Next we choose $i \not \equiv \pm j \mod q$, let 
\[ \psi_i = \phi_i \cdot \phi_{q-i} \quad \hbox{ and } \quad \psi_j = \phi_j 
\cdot \phi_{q-j} , \]  
and verify
\begin{eqnarray}
\label{todiv} 
\psi_i - \psi_j & = & \eta_{i,j} \cdot \lambda^{-2} \cdot
(\rho_1+\overline \rho_1)^2 \\ & = & \eta_{i,j} \cdot \lambda^{-2} \cdot
\phi_0^2 \cdot \left(\frac{ q^{nq-1} \cdot (\zeta^{1/2}+\overline
\zeta^{1/2})}{\varepsilon}\right)^2. \nonumber
\end{eqnarray}
We let $X = \mu_i \cdot \mu_{q-i}, Y = -\mu_j \cdot \mu_{q-j}$ and $Z =
\mu_0^2$. While $\mu_i$ are imaginary numbers, $X, Y$ are real, so $X,
Y \in \KL_p \cap \R = \KL^{++}$; trivially, $Z \in \Z[\zeta+\overline
\zeta] \subset \KL^{++}$. Now write $q = \varepsilon_1 \cdot
\lambda^{q-1}$ for the obvious real unit $\varepsilon_1 \in
\id{O}(\KL^{++})^{\times}$, set $N = 2 \left((q-1)(nq-1)-1\right)$ and 
\[ \delta = \eta_{i,j} \cdot \left(\frac{\zeta^{1/2}+\overline
\zeta^{1/2}}{\varepsilon}\right)^2 \cdot \varepsilon_1^N \in
\id{O}(\KL^{++})^{\times} . \] Inserting these new notations in
\rf{todiv} leads, after division by $\eta_0^2$, to:
\[ X^q + Y^q = \delta \cdot \lambda^N \cdot Z^q. \]
Once again we can apply Theorem \ref{tdesc}, obtaining a contradiction
with $q \nmid h_{pq}^-$. This completes the proof of this case.
\end{proof}
\subsection{The Ast\'{e}risque Case}
We denote the case $e=f=0$ by \textit{ast\'erisque} case. This is the
only case in which our results still depend on $x$ and $y$ -- the
obstruction to a more general result.

We suppose that $x \equiv 0 \mod q^{2}$ and since $e=0$, then $A^{q} =
y+x =_{q} y$, so that $y$ is a $q$-adic \nth{q} power. By \rf{sixqpow}
there is in this case a $\rho \in \Z[\zeta]$ such that $\rho^{q} = y
+ x \cdot \overline \zeta^{2q}$ and
\[ \frac{\zeta^{q} \cdot \rho^{q} + \overline{\zeta^{q} \cdot
    \rho^{q}}}{\zeta^{q} + \overline \zeta^{q}} = x+y = A^{q} . \]
By the usual argument of Lemma \ref{idealq}, we have then
\[ \phi_{i} = \frac{\zeta \xi^{i} \cdot \rho + \overline{\zeta
    \xi^{i} \cdot \rho}}{\zeta \xi^{i} + \overline {\zeta \xi^{i}}} =
\delta_{i} \cdot \mu_{i}^{q}, \] for some units $\delta_{i} \in
\Z[\zeta, \xi]$ and $i = 0, 1, \ldots, q-1$. Next, we investigate
these units $q$ - adically. Let $t^{q} \equiv y \mod q^{N}$, for some
integer $t \equiv A \mod q^{v_{q}(x)-1}$, or, likewise, $t$ be a $q$ -
adic approximation of the \nth{q} root of $y$. Then
\[ \delta_{i} =_{q} t \cdot \left(1 + \frac{\zeta^{1-2q} \cdot \xi^{i}
    + \overline{\zeta^{1-2q} \cdot \xi^{i}}}{\zeta \xi^{i} + \overline
    {\zeta \xi^{i}}} \cdot \frac{x}{q y}\right).\] In particular
\[ \delta_{0} =_{q} t \cdot \left(1 + \frac{\zeta^{2q-1}+\overline
    \zeta^{2q-1}}{\zeta+\overline\zeta } \cdot \frac{x}{q \cdot y}
\right) \in \Z[\zeta]. \] 

Let $\lambda = (1-\xi)$ (note the deviation from the usual definition
of $\lambda!$), so $\xi^{i} \equiv 1 -i \lambda \mod \lambda^{2}$. A
further investigation of the units $\psi_{i} = \delta_{i}/\delta_{0}$
shows that $\psi_{i} =_{q} 1 + \frac{x}{q y} \cdot c(\zeta) \cdot
\lambda + O(q \lambda^{2})$, with
\begin{eqnarray*}
 \lambda c(\zeta)  \equiv  \lambda \frac{\zeta^{1-2q} \cdot \xi^{i}
    + \overline{\zeta^{1-2q} \cdot \xi^{i}}}{\zeta \xi^{i} + \overline
    {\zeta \xi^{i}}} - \frac{\zeta^{2q-1}+\overline
    \zeta^{2q-1}}{\zeta+\overline\zeta } 
 \equiv  2 i \cdot \lambda \cdot \frac{\zeta^{2q}-\overline
    \zeta^{2q}}{(\zeta+\overline \zeta)^{2}} \mod \lambda^{2},  
\end{eqnarray*}
and thus $c(\zeta) = 2i \cdot \frac{\zeta^{2q}-\overline
  \zeta^{2q}}{(\zeta+\overline \zeta)^{2}} \neq 0$.  Thus Lemma
\ref{unitq1} implies that $\psi_{i}$ must be a \nth{q} power and
$c(\zeta) \equiv \sigma_{q}(\beta) - \beta \mod q$, for some $\beta
\in \Z[\zeta]$. In particular, if $q \equiv 1 \mod p$, then
$\sigma_{q}(\beta) \equiv \beta \mod \eu{Q}$ for all degree one primes
$\eu{Q} | q$ of $\Z[\zeta]$. But then $c(\zeta) \equiv 0 \mod q$,
which is in contradiction with $\Norm_{\K/\Q}(c(\zeta)) = (2i)^{p-1}
\cdot p \not \equiv 0 \mod q$. In this case we should have $x/q \equiv
0 \mod q^{2}$. A fortiori, $\delta_{0} \equiv t \mod q^{2}$ and so by
Proposition \ref{un1q} it must be a \nth{q} power. We have $\delta_{0}
\equiv t \equiv A \mod q^{2}$ and, since it is a \nth{q} power, also
$\delta_{0} = \gamma^{q}$ for some unit $\gamma$. But then $y \equiv
\gamma^{q^{2}} \mod q^{3}$, and so $y$ is a $q$ - adic \nth{q^{2}} power
with $y^{p-1} \equiv \Norm(\gamma^{q^{2}}) \equiv 1 \mod q^{3}$.   

One notes that if $-1 \in <q \mod p>$ and thus $q$ splits in real
primes in $\K$, then $c(\zeta)$ always verifies the condition $c =
\sigma_{q}(\beta) - \beta$ - since in fact $c(\zeta)+\overline
c(\zeta) = 0$ and the congruence holds modulo the real primes above
$q$. What can be said more generally, if $-1 \not \in < q \mod p> $?
One answer is that one can prove the statement of Lemma \ref{ast1} in
this case, provided that additionally $p \equiv 1 \mod 4$. We have thus
our first partial result:
\begin{lemma}
\label{ast1}
Suppose that \rf{FC} has a solution with $e = f = 0$ and $q \equiv 1
\mod p$ of $-1 \not \in <q \mod p>$ and $p \equiv 1 \mod 4$. Then
$q^{3} | x$ and $y^{p-1} \equiv 1 \mod q^{3}$.
\end{lemma}
\begin{proof}
  The case $q \equiv 1 \mod p$ was already explained above. If $-1
  \not \in < q \mod p>$, then let $g$ generate $\ZMs{p}$ and $H = \{
  g^{i} : 0 \leq i < (p-1)/2 \}$ be a set of represnetatives of
  $\ZMs{p}/\{-1,1\}$, and $H' = \ZMs{p} \setminus H$. Thus $x \in H
  \Leftrightarrow (p-x) \in H'$ and by hypothesis, $<q \mod p> \subset
  H$. Let $\eu{q}$ be some prime above $p$, so $\eu{q} \neq \overline
  {\eu{q}}$. The condition $c(\zeta) = \sigma_{q}(b)-b$ implies then
  $\sum_{x \in H} \sigma_{x}(c(\zeta)) \equiv \sum_{x \in < q \mod p>}
  \sigma_{x}(c(\zeta) \equiv 0 \mod q$. Let 
\[ a =
\frac{\zeta^{2q}}{(\zeta+\overline \zeta)^{2}} =
\frac{\zeta^{2(q+1)}}{(1+\zeta^{2})^{2}} = \sigma_{2}
\left(\frac{\zeta^{q+1}}{(1+\zeta)^{2}}\right) ,\]
so that $c(\zeta) = a - \overline a$ and the previous condition
amounts to 
\begin{eqnarray}
\label{con1}
\sum_{x \in H} \sigma_{x}(a) \equiv \sum_{x \in H'} \sigma_{x}(a) \mod
q. 
\end{eqnarray}
We shall show that this condition cannot be fulfilled if $p \equiv 1
\mod 4$.
\end{proof}
\begin{remark}
  The above result also implies that for all $q$ we have $\phi_{i} =
  \delta_{0} \cdot \mu_{i}^{q}$ and $\mu_{0}^{q} =
  \mu_{i}^{q}+\mu_{q-i}^{q}$; this fact is noteworthy but leads
  unfortunately to no further descent.
  
  It may also be observed that if $-1 \in <q \mod p>$ and thus $q$
  splits in real primes in $\K$, then $c(\zeta)$ always verifies the
  condition $c = \sigma_{q}(\beta) - \beta$. We have already shown
  that this is not the case if $q \equiv 1 \mod p$. What can be said
  more generally, if $-1 \not \in < q \mod p> $? One answer is that
  one can prove the statement of Lemma \ref{ast1} in this case,
  provided that additionally $p \equiv 1 \mod 4$.
\end{remark}
We proceed with some global estimates. We shall assume that $|x| >
|y|$, which is allowed since $x, y$ are interchangeable for the global
estimates; furthermore, we assume that $x > 0$. Note that this choice
does not allow any more to choose which of $x$ and $y$ is divisible by
$q^{3}$; this is of no relevance for the global estimates we are about
to prove. We have
\begin{lemma}
\label{astyest}
Suppose that \rf{FC} has a solution with $e = f = 0$ and $x > |y| >
0$. Then 
\begin{eqnarray}
\label{yest}
|y| > c(q) \cdot x^{1-2/q}, 
\end{eqnarray}
for some absolutely computable, strictly increasing function $c(q)$
with $c(5) > 1$.
\end{lemma}
\begin{proof}
  Suppose first that $y > 0$ and let $\psi = \rho \cdot \overline \rho
  = (x+\zeta y)(x+\overline \zeta y) \in \R \cap \Z[\zeta]$. Then
\[ \psi^{q} = \left((x+y)^{2} - \mu x y\right) = A^{2q} \cdot
\left(1 - \frac{ \mu \cdot x y}{(x+y)^{2}}\right) \quad \hbox{
  with } \quad \mu = (1-\zeta)(1-\overline \zeta) .\] Note that
$\left| \frac{ \mu \cdot x y}{(x+y)^{2}}\right| \leq |\mu|/4
< 1$ in this case, so there is a converging global binomial expansion
of $f(x, y) = \left(1 - \frac{ \mu \cdot x
    y}{(x+y)^{2}}\right)^{1/q}$. The expressions $\psi$ and $A^{2}
  f(x, y)$ have the same \nth{q} power, so they differ by a \nth{q}
  root of unity. But since they are both real, they must coincide:
  $\psi = A^{2} \cdot f(x, y)$.  Furthermore, the
series summation commutes with the action of $\Gal(\Q(\zeta)/\Q)$, for
the same reason \footnote{ \ For more detail on this kind of argument,
  see for instance \cite{Mi}}. Thus, for all $\sigma \in G_{p}$,
\[ \sigma(\psi) = A^{2} \cdot \left(1 + \sum_{n=1}^{\infty}
  \binom{1/q}{n} \cdot \left(\frac{-\sigma(\mu) \cdot x
      y}{(x+y)^{2}}\right)^{n} \right) . \] An easy computation (see
e.g. \cite{Mi}) shows that the binomial coefficients are bounded by
$\left|\binom{1/q}{n} \right| < \frac{1}{qn}$, while
$|\sigma(\mu)| < 4$ for all $\sigma \in \Gal(\Q(\zeta)/\Q)$. If
$R_{2}$ is the second order remainder of the above series, one finds
from these estimates that
\[ \left|\sigma R_{2} \right| < \left(\frac{4 A \cdot x
  y}{(x+y)^{2}}\right)^{2} \cdot \frac{2}{q}
  \ln\left(\frac{x+y}{x-y}\right), \]
uniformly for $\sigma \in \Gal(\Q(\zeta)/\Q)$.

For a fixed $\sigma_{0} \in \Gal(\Q(\zeta)/\Q)$, we now give a uniform
estimate of the difference $\delta = |\psi - \sigma_{0}(\psi)| \in
\Z[\zeta]$. Since $\delta | \left(\psi^{q} - \sigma_{0}(\psi)^{q}
\right) = \left(\mu-\sigma_{0}(\mu)\right) \cdot x y \neq
0$, it follows that $\delta$ is a non vanishing algebraic integer. Its
absolute value is:
\begin{eqnarray*} 
|\delta| & = & \left| A^{2} \cdot \frac{(\mu-\sigma_{0}(\mu))
    x y}{q (x+y)^{2}} + (R_{2} - \sigma_{0}(R_{2}))\right| \\
& < & A^{2} \cdot \frac{4xy}{q(x+y)^{2}} + 2 \left(\frac{4 A \cdot x
  y}{(x+y)^{2}}\right)^{2} \cdot \frac{2}{q}
  \ln\left(\frac{x+y}{x-y}\right) \\
& = &\frac{4 A^{2} x y}{q(x+y)^{2}} \cdot \left(1 + \frac{16
    xy}{(x+y)^{2}} \cdot \ln\left(\frac{x+y}{x-y}\right)\right).
\end{eqnarray*}
Note that the above estimate holds for all $\sigma \delta$ uniformly;
one then verifies that for $|y| < c(q) \cdot |x|^{1-2/q}$ and, say,
$c^{-1}(q) = 4/q \cdot \left(1 + 16 q^{-2} \cdot \ln\left(1+2/q^{q}
\right) \right)$ (use the lower bound on $|x|$, above!), then $0 <
|\sigma \delta| < 1$ and thus $\Norm |\delta| < 1$, in contradiction
with the fact that $\delta$ is a non vanishing algebraic integer. This
completes the proof for $y > 0$.

If $y < 0$, one lets $y' = -y$ in the previous proof and sets
$\mu = (1+\zeta)(1+\overline \zeta)$. Concretely, we have
$\psi^{q} = (x+y')^{2} \cdot \left(1 - \mu
  \frac{xy}{x+y'}\right)$. As a result, the factor $(x+y')^{2}$ is not
the \nth{q} power of an integer any more, one replaces $A^{2}$ by the
real positive value of $(x-y)^{2/q}$. This has no impact on the
estimates of the algebraic integer $\delta$, which are perfectly
analog, and lead to the same result.
\end{proof}
\subsection{The Proof of Theorem \ref{main2}}
\begin{proof}
  Suppose that \rf{FC} has a non trivial solution for odd primes $p,
  q$ verifying the conditions of the Theorem. Then by Theorem
  \ref{main} it follows that $p^{e} | z$ and $x+f y \equiv 0 \mod
  q^{2}$, for some $f \in \{-1,0,1\}$. The cases $f = 1$ are
  impossible, as proved in Theorems \ref{desc1} and \ref{desc2}. Three
  of the remaining cases are dealt by some Wieferich condition, as
  proved in Propositions \ref{p1}, \ref{p2} and \ref{p3}, while for
  the Ast\'erisque case we have shown in Lemma \ref{ast1} that $q^{3}|
  x$ if $q \equiv 1 \mod p$. We still have to prove the lower bound
  \rf{bound}.
 
  Since the lower bound on $\max\{|x|,|y|\}$ is slightly better in the
  particular case Ast\'erisque then the general bounds in Theorem
  \ref{main1}, we deduce this bound separately.
  
  We assume now that $q^{2}|x$ and $k = v_{q}(x) \geq 2$ - thus
  dropping the assumption $|x| > |y|$. By letting $t^{q} = y$ as
  elements of $\Z_{q}$, the $q$ - adic expansion of $\rho = (y + \zeta
  x)^{1/q}$ is then
\[ \rho = t \cdot \left(1 + \sum_{n=1}^{\infty} \binom{1/q}{n} \cdot
  (\zeta x/y)^{n} \right) . \] If $A^{q} = x+y$, we now consider the
algebraic integer $\delta = A + \Tr_{\Q(\zeta)/\Q} ( \zeta \rho) \in
\Z$. Since $A^{q} = y(1+x/y)$, one observes that the $q$-adic expansion
of $A$ results from the one for $\rho$ by replacing $\zeta$ with $1$.
An easy computation yields (note the factor $\zeta$ of $\rho$ in the
definition of $\delta$) the $q$ - adic expansion:
\[ \delta = p t \cdot \binom{1/q}{p-1} \cdot (x/y)^{p-1} +
O\left(\binom{1/q}{p} \cdot (x/y)^{p}\right) .\] Obviously, $\delta
\neq 0$. In order to see this, note that $V(j) =
-v_{q}\left(\binom{1/q}{j}\right) = j+v_{q}(j!)$ by \rf{binomvq}.
Furthermore, let $W(j) = v_{q}\left(\binom{1/q}{j} \cdot (x/y)^{j}
\right) = k j - V(j)$ and thus
\[ W(p) - W(p-1) = k + (V(p-1) - V(p)) = k-1+v_{q}((p-1)!)-v_{q}(p!) =
k-1 > 0 . \] It follows that $\delta \equiv \binom{1/q}{p-1} \cdot
(x/y)^{p-1} \mod q^{W(p)}$ and thus $\delta \neq 0$. It follows in
particular that $\delta \equiv 0 \mod q^{(k-1)(p-1)-v_{q}((p-1)!)}$. Let
now $B = |x|+|y|$; then $A^{q} \leq B$ and $|\sigma(\zeta \rho)|^{q} < B$, so
$|\delta| < p B^{1/q}$. The last inequalities combine to:
\begin{eqnarray}
\label{bd1}
 \max(|x|, |y|) > (|x|+|y|)/2 > \frac{1}{2} \cdot
\left(\frac{q^{(k-1)(p-1)}}{p} \right)^{q} .
\end{eqnarray}
If $|x| > |y|$, then one can use the bound in \rf{bd1} which is
stronger then \rf{bound}. Otherwise, \rf{bd1} implies $|y| >
\frac{1}{2} \cdot \left(\frac{q^{(k-1)(p-1)}}{p} \right)^{q}$ and by
interchanging $x$ and $y$ in \rf{yest} (the maximum is now $|y|$), we
obtain the claim \rf{bound}, where $c_{1}(q) = c(q)/2$, with $c(q)$
from \rf{yest}. If $q \not \equiv 1 \mod p$, all we know is $k \geq
2$, which yields \rf{bound}; otherwise, by Lemma \ref{ast1}, we have
$k \geq 3$ and \rf{bound1}.

This result improves upon \rf{lowb}. It is however due to \rf{yest}
that one obtains a lower bound on $\min\{|x|, |y|\}$, which allows to
assert that $x^{p}+C^{p} = z^{q}$ has no solutions for $|C|$ below
this lower bound, as we have explicitly shown in the Corollary
\ref{catg1}.
\end{proof}

\section{The Equation of Catalan in the Rationals}
We have proved in Lemma \ref{rcfc} that the rational Catalan equation
\rf{ratcat} is equivalent to \rf{homfc}:
\[ X^{p}+Y^{q} = Z^{pq} .\]
Note that this equation is \textit{symmetric} in $p, q$ in the sense
that it splits in the two equations:
\begin{eqnarray}
\label{eqp}
X^{p} + (-Z^{q})^{p} & = & (-Y)^{q}, \\
\label{epq}
Y^{q} + (-Z^{p})^{q} & = & (-X)^{p}, 
\end{eqnarray}
which are both of type \rf{FC}. Thus Theorem \ref{main2} applies to
both equations. The task we still have to achieve for proving Theorem
\ref{trc} consists in eliminating the Ast\'erisque Case, by using the
symmetry in the above equations. This is a consequence of the
following:
\begin{proposition}
\label{AstCat}
Suppose that $p, q$ are two odd primes verifying the premises of
Theorem \ref{main2} and for which \rf{eqp} holds. Then $q | X$.
\end{proposition}
\begin{proof}
  Under the given premises, Theorem \ref{main2} implies that $q |
  \left(X \cdot Z^{q}\right)$. For clarity, we use the substitution $x
  = X, y = (-Z)^{q}$ and $z = -Y$ in order to bring \rf{eqp} in the
  form of the reference Fermat - Catalan equation \rf{FC}. We will
  show that the assumption $q | Z$ -- and thus $q | y$ -- leads to a
  contradiction.  For this we use again the Descent Theorem and the
  fact $y = (-Z)^{q}$ is a \nth{q} power.

We assume thus that $q | y = (-Z)^{q}$ and $p \nmid xyz$. From \rf{sixqpow}
we have in this case $\rho^{q} = x + \zeta y$ and thus
\[ (\zeta - \overline \zeta) y = -(\zeta - \overline \zeta) Z^{q} =
\rho^{q}-\overline \rho^{q}. \] Let $\phi'_{i} = \xi^{i} \rho -
\overline \xi^{i} \overline \rho$ and $\phi'_{0} = \rho-\overline
\rho$. Then $\prod_{i=0}^{q-1} \phi'_{i} = (\zeta - \overline \zeta)
y$ and since $(y, p) = 1$ and $\phi'_{i} = \tau_{i}(\phi'_{1})$, while
$p \not \equiv 1 \mod q$ and thus $\wp = (1-\zeta)$ does not split
completely in $\KL/\K$.  It follows that $\wp | \phi'_{0}$ and
$(\phi'_{i}, \wp) = (1)$. Let $y = q^{nq} \cdot C^{q}$, with $(C, pq) =
1$. By introducing the normalization
\[ \phi_{0} = \frac{\rho-\overline \rho}{q^{nq-1} \cdot (\zeta
  -\overline \zeta)} \quad \hbox{ and } \quad \phi_{i} = \frac{\xi
  \rho - \overline \xi \overline \rho}{\xi - \overline \xi}, \] the
arguments use in the proof of the Descent Theorem yield here:
\begin{eqnarray*} 
\prod_{i=0}^{q-1} \phi_{i} & = & C^{q},  \quad \hbox{and} \\  
(\phi_{i}, p \cdot q) & = & (\phi_{i}, \phi_{j}) = (1), \quad i \neq
j \geq 0 .
\end{eqnarray*}  
We can apply now Lemma 5 and find that $\phi_{i} = \delta_{i} \cdot
\mu_{i}^{q}$, for $i = 0, 1, \ldots, q-1$. If $t \in \Z_{q}$ is such
that $t^{q} = x$ (existence is provided by Proposition \ref{p1}), then
the $q$ - adic expansion of $\delta_{i}$, given that $q^{2} | y$,
yields: $\delta_{i} =_{q} t$. This must then be a \nth{q} power, by
Lemma \ref{unitp} and it follows plainly that $\phi_{i} =
\mu_{i}^{q}$.  The proof proceeds like in the one for the first
descent case and shall be sketched here. We define $\psi_{i} =
\phi_{i} \cdot \phi_{q-i} = (\rho+\overline
\rho)^{2}/(\xi^{i}-\overline \xi^{i})^{2} - \rho \cdot \overline \rho
$ and find that
\[ \psi_{i} - \psi_{j} = (\mu_{i} \cdot \mu_{q-i})^{q} - (\mu_{j}
\cdot \mu_{q-j})^{q} = \delta_{i,j} (\zeta - \overline \zeta)^{2}
\cdot q^{2(nq-1)} \cdot \mu_{0}^{2q} ,\] with $\delta_{i,j}$ defined
in the proof of Theorem \ref{desc1}, so that $(\xi-\overline \xi)^{2}
\delta_{i, j} \in \Z^{\times}[\zeta, \xi]$. The descent argument is in
place and the claim of our Theorem follows from the assumption by
means of Theorem \ref{desc}.
\end{proof}
\subsection{Proof of Theorem \ref{trc}}
We can now complete the proof of Theorem \ref{trc}.
\begin{proof}
We know by Lemma \ref{rcfc} that \rf{eqp} and \rf{epq} hold
simultaneously. The additional conditions ensure that the premises of
Theorem \ref{main2} hold for both equations, considered as equations
of the type \rf{FC} (e.g. by substitutions like in the proof of the
previous Proposition). We analyze the consequences of the six
conditions in Theorem \ref{trc}; for this we refer the reader to the
case analysis made for the proof of Theorem \ref{main2}.

The conditions 1., 2. and 6. are sufficient for eliminating the
descent cases $f = 1$ in both equations \rf{eqp} and \rf{epq}. The
conditions 3. and 4. then show that the cases with $f = -1$ cannot
occur for either \rf{eqp} or \rf{epq}. The only cases left are thus
the ones with $f = 0$. Finally, condition 5. implies that the case $e
= 1, f = 0$ does also not occur and the only case left is the
Ast\'erisque case $e = f= 0$, for both \rf{epq} and \rf{eqp}. However,
by Proposition \ref{AstCat}, this implies that $q | X$. But this is
exactly the case $e = 1$ in \rf{epq}, which is granted not to have
solutions by the same condition 5. The contradiction completes the
proof of the Theorem. 
\end{proof}

\vspace*{0.3cm}

%{\bf Acknowledgments:}
%  \vspace*{0.5cm}

%

\begin{thebibliography}{ZZZZZZ}

\bibitem[Be]{Be} Frits Beukers: {\em The diophantine equation
    $Ax^p+By^q=Cz^r$}, Lectures held at Institut Henri Poincare,
    September 2004, http://www.math.uu.nl/people/beukers/

\bibitem[Br]{Br} Niels Bruin: {\em The Diophantine equations $x^{2}
    \pm y^{4} = \pm z^{6}$ and $x^{2}+y^{8}=z^{3}$}, Compositio
    Math. \textbf{118} (1999) no. 3, pp. 305-321.

\bibitem[Br2]{Br2} Niels Bruin: {\em Chabauty Methods using Elliptic
    Curves}, J. f\"ur die reine und angewandte Mathematik, {\bf 562}
    (2003), pp. 27-49.

\bibitem[BH]{BH} Y. Bugeaud and G. Hanrot: {\em Un nouveau crit\`ere pour
    l'\'equation de Catalan}, (2000), Matematika {\bf 47} (2000), pp.
  15-33.

\bibitem[Bu]{Bu} Y. Bugeaud: Private communication.
  
\bibitem[Da]{Da} Henri Darmon: {\em Rigid local systems, Hilbert
    modular forms, and Fermat's Last Theorem}. Duke Math J.
  \textbf{102} (2000) 413-449.
  
\bibitem[Da1]{Da1} Henri Darmon: {\em The equations
    $x^{n}+y^{n}=z^{2}$ and $x^{n}+y^{n}=z^{3}$}, Internat.  Math.
  Res. Notices 1993, \textbf{10} 263--274.
  
\bibitem[Da2]{Da2} Henri Darmon: {\em Modularity of fibers in rigid
    local systems}, Annals of Math, \textbf{149} (1999) 1079-1086.

\bibitem[DG]{DG} H. Darmon and A. Granville: {\em On the Equation
    $z^{m} = F(x,y)$ and $Ax^{p}+By^{q}=Cz^{r}$}, Bull. London
    Math. Society, \textbf{27} (1995), no. 6, pp. 513- 543.
    
  \bibitem[DM]{DM} Henri Darmon and Loic Merel: {\em Winding quotients
      and some variants of Fermat's last theorem}. J. Reine Angew.
    Math. \textbf{490} (1997), 81--100.  
    
    \bibitem[El]{El} Jordan Ellenberg: {\em Galois representations
      attached to $\Q$-curves and the generalized Fermat equation
      $A^{4} + B^{2} = C^{p}$}, Amer.  J.  Math. \textbf{126, (4)},
    763--787 (2004).

\bibitem[Kr]{Kr} A. Kraus: {\em Sur l'\'equation $a^{3}+b^{3} =
    c^{p}$} Experimental Math. {\bf 7} (1998), pp. 1-13.
 
\bibitem[La]{La} Lang, S.: {\em Algebraic Number Theory}, 
Second Edition, Springer 1986, Graduate Texts in 
Mathematics \textbf{ 110 } 

\bibitem[Lo]{Lo} R. Long: {\em Algebraic number theory}, Marcel
  Dekker, Monographs in Pure and Applied Mathematics {\bf 41}, (1977).

\bibitem[Mi]{Mi} P. Mih\u{a}ilescu: {\em Primary Cyclotomic Units
and a Proof of Catalan's Conjecture}, Crelle's Journal, \textbf{572}
(2004), pp. 167-195.

\bibitem[Mi1]{Mi1} P. Mih\u{a}ilescu: {\em A Cyclotomic Approach to
    the Fermat -- Catalan Conjecture}, Octogon (mathematical magazine,
  M. Bencze ed., Bra\c{s}ov, Rom\^{a}nia), \textbf{12}, No. 1, (2004),
  p. 5-11
  
\bibitem[Mi2]{Mi2} P. Mih\u{a}ilescu: {\em Reflection, Bernoulli
    numbers and the proof of Catalan's conjecture}, to appear in
  Proceedings of the Fourth European Congress of Mathematics,
  Stockholm 2004.
  
\bibitem[PSS]{PSS} B. Poonen, E. Schaefer and M. Stoll: {\em On
    $x^{2}+y^{3}+z^{7} = 0$}, Presentation at a Seminary in
  Oberwolfach, to appear.

\bibitem[Ri]{Ri} P. Ribenboim: {\em $13$ Lectures on Fermat's Last
Theorem}, Springer Verlag (1979).

\bibitem[ST]{ST} T.N. Shorey and R. Tijdeman: {\em Exponential
    diophantine equations}, Cambridge University Press, 1986.

\bibitem[Wa]{Wa} L. Washington: {\em Introduction to 
Cyclotomic Fields}, Second Edition, Springer (1996), 
Graduate Texts in Mathematics {\bf 83}.

\bibitem[Za]{Za} D. Zagier, personal communication.
\end{thebibliography}
\end{document}